\documentclass[12pt]{amsart}
\usepackage{a4,amsthm, epsf, amscd}

\theoremstyle{plain}
\newtheorem{theorem}{Theorem}[section]
\newtheorem{lem}[theorem]{Lemma}
\newtheorem{rem}[theorem]{Remark}
\newtheorem{prop}[theorem]{Proposition}
\newtheorem{cor}[theorem]{Corollary}
\newtheorem{defi}[theorem]{Definition}
\theoremstyle{definition}

\def\Exp{{\rm Exp}}
\def\ad{{\rm ad}}
\def\Ad{{\rm Ad}}
\def\a{{\mathfrak{a}}}
\def\m{{\mathfrak{m}}}
\def\u{{\mathfrak{u}}}
\def\p{{\mathfrak{p}}}    
\def\k{{\mathfrak{k}}}
\def\g{{\mathfrak{g}}}
\def\h{{\mathfrak{h}}}

\def\C{{\mathbb{C}}}
\def\R{{\mathbb{R}}}

\def\nsmallskip{\smallskip\noindent}
\def\bbigskip{\bigskip\bigskip}
\def\nbigskip{\bigskip\noindent}

\def\nmedskip{\medskip\noindent}
\def\buildunder#1#2{\mathrel{\mathop{\kern0pt #2}
\limits_{#1}}}
\def\buildover#1#2{\buildrel#1\over#2}

\def\qq{/\kern-.185em /}

\begin{document}

\title[A family of complexifications]{A family of adapted complexifications
for $\,SL_2(\R)$}
 
\author[Halverscheid]{S. Halverscheid*}
\author[Iannuzzi]{A. Iannuzzi}

\address{Stefan Halverscheid: Universit\"at Bremen
Bibliothekstr. 1,
D-28359 Bremen, Germany}
\email{sth@math.uni-bremen.de}

\address{Andrea Iannuzzi: Dip. di Matematica,
II Universit\`a di Roma  ``Tor Vergata", Via della Ricerca Scientifica,
I-00133 Roma, Italy}
\email{iannuzzi@mat.uniroma2.it}

\thanks{*\,\,Partially supported by DFG-Forschungsschwerpunkt 
``Globale Methoden in der komplexen Geometrie".}

\begin{abstract}
Let $\,G\,$ be a non-compact, real semisimple Lie group.
We consider maximal complexifications of $\,G\,$ which are adapted to a 
distinguished one-parameter family of naturally reductive,
left-invariant metrics. In the case of $\,G=SL_2(\R)\,$
their realization as equivariant Riemann domains over $\,G^\C=SL_2(\C)\,$
is carried out and their complex-geometric properties are investigated.
One obtains new examples of non-univalent, non-Stein,
maximal adapted complexifications.
\end{abstract}


\maketitle

\section{Introduction}

\medskip
\smallskip
Let $\,\nabla\,$ be a linear connection on
a real-analytic manifold $\, M\,$ which is identified with the zero section in the tangent bundle $\,TM$.
A complex structure defined on a domain $\,\Omega\,$ of $\,TM\,$ containing $\,M\,$
 is adapted to  the connection if 
for  any $\,\nabla$-geodesic $\,\gamma\,$
its complexification $\,\gamma_*$, given by
$\, (x + iy) \to y \, \gamma'(x)$,  
is holomorphic on  $\,(\gamma_*)^{-1}(\Omega)$. In this situation we refer 
to $\,\Omega\,$ as an adapted complexification of $\,(M,\nabla)$. 
In the case where $\,\nabla\,$ is real-analytic, 
R. Bielawski (\cite{Bi}) and R. Sz{\H o}ke  (\cite{Sz3})
recently showed that the adapted complex structure exists in a
neighborhood of $\,M$.
In \cite{Sz3} one also finds a uniqueness result.
For the Levi Civita connection of a real-analytic Riemannian 
manifolds, such results were known since the pioneering works of
Guillemin-Stenzel (\cite{GuSt}) and Lempert-Sz{\H o}ke (\cite{LeSz}).

In the presence of a  ``large enough" Lie group acting on $\,M\,$ and
preserving the geodesic flow induced by $\,\nabla$, one can prove that there
exists a maximal domain
$\,\Omega\,$ for the adapted complex structure, i.e.$\,$every 
adapted complexification is necessarily contained in $\,\Omega\,$
(see$\,$Prop.$\,$3.1, cf. \cite{HaIa1}).
If $\,M\,$ is a non-compact, Riemannian symmetric space such a maximal complexification is
well-known under the name of Akhiezer-Gindikin domain (see \cite{AkGi}, cf. \cite{FHW}).

Recall that $\,M\,$ is the fixed point set of the anti-holomorphic involution
on $\,\Omega \subset TM\,$
given by $\,v \to -v\,$ and in the case of the Levi Civita connection 
associated to a pseudo-Riemannian manifold the metric $\,\nu\,$ appears as the 
restriction of a pseudo-K\"ahler metric $\,\kappa\,$ with the same 
index as $\,\nu$.
Moreover $\,\kappa\,$ admits a global potential
whose properties give important geometric insights of
$\,\Omega\,$ ([Sz3], cf. \cite{Bu}, \cite{LeSz}, \cite{PaWo}, \cite{St}, \cite{Sz1}).

For a  connected, non-compact, real semisimple Lie group $\,G$, let
$\,\g=\k \oplus \p\,$ be the Cartan decomposition of its 
Lie algebra $\,\g\,$ with respect to a maximal compact subalgebra 
$\,\k$. Denote by 
$\,B\,$  the Killing form of $\,G\,$ and consider 
the distinguished one-parameter family of left-invariant metrics
$\,\nu_m\,$ (\,degenerate for $\,m=0\,$) uniquely defined by
$$\,\nu_m|_{\g} (X,Y)= -mB(X_\k,Y_\k) + B(X_\p,Y_\p)\,,$$

\nsmallskip
for any  
$\,X=X_\k+X_\p\,$ and $\,Y=Y_\k+Y_\p\,$ in $\,\k \oplus \p$.
Note that for all real $\,m$,  the action of  the group $\,L:=G \times K$,
given by left and right multiplication on $\,G$, is isometric.
Here $\,K\,$ is the connected subgroup of $\,G\,$ generated by $\,\k$.
Recall that for $\,m>0\,$ these metrics appear
in the classification of all naturally 
reductive,  left-invariant Riemannian metrics on $\,G\,$ given by C. Gordon in \cite{Go}.

Our main goal is to present 
new examples of maximal complexifications adapted, in the non-degenerate cases, to
the Levi Civita connection associated to a metric of the above family. In the degenerate 
case one finds a maximal complexification adapted to the unique real-analytic linear connection which is obtained as the limit of such Levi Civita connections. 
For $\,G=SL_2(\R)$
we give a precise description of these complexifications
and we determine their 
basic complex-geometric
properties.
For positive $\,m\,$ this gives, 
along with previous results (see \cite{Sz2}, \cite{BHH}, \cite{HaIa2}),
examples among all classes of $\,3$-dimensional,
naturally reductive, Riemannian homogeneous spaces (cf. \cite{BTV}).
For  $\,m=-1\,$ one obtains the symmetric pseudo-Riemannian case which
has been investigated, among others, by G. Fels (\cite{Fe})
and R. Bremigan  (\cite{Br}).

The paper is organized as follows.
Basic results and properties of the above
metrics $\,\nu_m\,$ are recalled in section 2. There we also point out 
with an example that in order to perform 
(pseudo) K\"ahlerian reduction in this
pseudo-Riemannian context, one may need more 
conditions than those necessary in the Riemannian case 
(Rem. \ref{AGUILAR}, cf. \cite{Ag}).

In section 3 we give a version of the characterization 
of maximal adapted complexifications given in
\cite{HaIa1} which is suitable in our situation (Prop.$\,$\ref{POLARMAP}). In particular we show the existence
of 
 the maximal  complexification $\,\Omega_m\,$ adapted
 to $\,\nu_m$. This is realized as 
an $\,L$-equivariant  Riemann domain 
over the universal complexification $\,G^\C\,$ 
of $\,G$,
with  {\it polar map} $\,P_m:\Omega_m \to G^\C\,$.
By considering the 
usual identification $\,TG \cong G \times \g$, such complexification can be 
described via a slice for the induced  
$\,L$-action by
$$\Omega_m \ = \ L \cdot  \Sigma_{m} \,,$$

\noindent
where $\,\Sigma_{m}  \subset\{e\}\times \g\,$
is  a semi-analytic subset of the product of $\,\k\,$ and 
the closure of a Weyl chamber in a maximal abelian subalgebra 
$\,\a\,$ of $\,\p$.  That is, $\,(e,X) \in \Sigma_{m}\,$ if and only if $\,X\,$
belongs to the intersection of sublevel sets of 
certain real-analytic functions of $\,\k \oplus \a\,$
(Prop. \ref{SLICE} and \ref{DOMAIN}).

The case $\,G=SL_2(\R)\,$ is carried out in detail and
the defining functions for $\,\Sigma_{m}$ as well as the polar map
 $\,P_m\,$  are explicitly
determined in terms of a fixed basis of 
$\,\k \oplus \a\,$ in sections 4 and 6. For this it is useful to have
concrete realizations of slices and quotients of 
$\,SL_2(\C)\,$ with respect to the involved actions, i.e.
those of $\,SL_2(\R)$, $\,SL_2(\R) \times SO_2(\R)\,$ and
$\,SL_2(\R) \times SO_2(\C)$.  
This is separately discussed in section 5.

Finally, in  sections 7 and  8  we single out the following different situations, 
whose boundary cases are given by the symmetric pseudo-Riemannian
$\,m=-1\,$ and the degenerate $\,m=0$.
For $\,m <-1\,$ all maximal adapted complexifications
are biholomorphic via $\,P_m\,$ to a non-Stein $\,L$-invariant domain.
Namely $\,SL_2(\C)\,$ without a single $\,SL_2(\R) \times SO_2(\C)$-orbit
(Thm.$\,$\ref{DOMAIN<-1}).
For $\,-1\le m \le 0\,$ the polar map $\,P_m\,$ remains injective but its
 non-Stein image misses more and more $\,SL_2(\R) \times SO_2(\C)$-orbits.
If
$\,m > 0$, the maximal adapted complexifications turn out to be neither
holomorphically convex, nor holomorphically separable (Thm.$\,$\ref{DOMAIN>-1}). 
In all cases the envelope of holomorphy of $\,\Omega_{m}\,$
is shown to be biholomorphic to $\,SL_2(\C)\,$ (cf. \cite{Sz3}, Sect. $\,$9).

Note that in the Riemannian context $\,m>0\,$ all metrics
$\,\nu_m\,$ have mixed
sign sectional curvature. A similar situation can be noticed 
in the examples discussed in  \cite{HaIa1}. In these examples certain left-invariant
Riemannian metrics on the generalized Heisenberg group were considered.
Their sectional curvature has mixed sign and the associated
maximal adapted complexifications have
similar complex-geometric properties.
It would be interesting to know if this is only a coincidence.

\nbigskip
{\it Acknowledgment.} We wish to thank R\'obert  Sz{\H o}ke for sharing
with us a preliminary version of  \cite{Sz3} and the referee for his helpful
suggestions.

\nbigskip


\section{Preliminaries}

\bigskip

Let $\,G\,$ be a non-compact, semisimple Lie group.
Here we recall basic properties of 
the one-parameter family of 
left-invariant  metrics on $\,G\,$ which will be considered in the sequel.
Such a family contains degenerate,  
Riemannian and  pseudo-Riemannian, naturally reductive 
metrics.
The Riemannian ones appear in the classification 
given by  C. Gordon in \cite{Go}.
More details and curvature computations
 can be found in \cite{HaIa2}.

We also recall those facts  
on adapted complex structures
which are 
needed in the present paper. Finally we
point out an example showing that 
the reduction procedure  indicated by R. Aguilar (\cite{Ag})
in the Riemannian context does not apply automatically when
dealing with pseudo-Riemannian geometry.


\bigskip
\begin{defi}
\label{NATRED}
 {\rm (cf. \cite{O'N})}  A pseudo-Riemannian metric $\,\nu\,$ on a 
 homogeneous manifold $\,M\,$ 
is naturally reductive if
there exist a connected Lie subgroup $\,L\,$ of $\,\mathrm{Iso}(M)\,$
acting transitively on $\,M\,$ and a decomposition 
$\,\mathfrak l =\mathfrak h \oplus \mathfrak m\,$ of $\,\mathfrak l$,
where $\,\mathfrak h\,$ is the Lie algebra of the isotropy
group $\,H\,$ at some point of $\,M$, such that 
$\,\Ad(H) \, \mathfrak m \subset \mathfrak m\,$
and 
$$\tilde \nu([X,\,Y]_\mathfrak m , \,Z) =
\tilde \nu(X,\, [Y,\,Z]_\mathfrak m)$$

\nsmallskip
for all $\,X,\,Y,\,Z \in \mathfrak m$. Here
 $\, [\ ,\ ]_\mathfrak{m}\,$ denotes the
$\mathfrak{m}$-component of $\, [\ ,\ ]\,$ and 
$\,\tilde \nu\,$ is the pull-back of $\,\nu\,$ to 
$\,\m\,$ via the natural projection $\,L \to L/H \cong M$.
In this setting we refer to $\,\mathfrak h \oplus \mathfrak m\,$ as a naturally reductive 
decomposition and to $\,L/H\,$ as a naturally reductive realization of $\,M$.
\end{defi}

\bigskip
For a naturally reductive realization $\,L/H\,$
every geodesic through the base point $\,eH\,$ 
is the orbit of a
one-parameter subgroup of $\,L\,$ generated by some
$\,X \,\in \,\mathfrak m\,$
(see \cite{O'N}, p. 313). In fact for a Riemannian homogeneous
manifold  $\,L/H\,$ with an $\,\Ad(H)$-invariant 
decomposition $\,\mathfrak h \oplus \mathfrak m\,$
this property implies that  $\,L/H\,$ is a
naturally reductive realization (see, e.g. \cite{BTV}).

\bigskip
Let $\,G\,$ be a connected, non-compact,  semisimple Lie group 
and  let $\,\g = \k \oplus \p\,$  be the Cartan decomposition
of its Lie algebra with respect to a maximal compact
Lie subalgebra $\,\k$.
Let $\,B\,$ denote the Killing form  on $\,\g\,$ and, for every 
real $\,m$, assign a left-invariant
metric $\,\nu_m\,$ on $\,G\,$ by defining its restriction on
$\,\g \cong T_eG\,$ as follows:  
\begin{equation}
\label{DEFMETR} 
\nu_m \bigr|_{\g} (X,Y)=
-mB(X_\k,Y_\k) + B(X_\p,Y_\p),
\end{equation}

\medskip
\noindent
for any $\,X=X_\k+X_\p\,$ and  $\, Y=Y_\k+Y_\p\,$ in  $\,\k \oplus \p$. 
Since $\,B\,$ is negative definite on $\,\k\,$ and positive definite on $\,\p$,
these metrics are Riemannian, degenerate or pseudo-Riemannian
when $\,m>0$, $\,m=0\,$ or  $\,m<0$,  respectively.

Let $\,K\,$ be the connected subgroup of $\,G\,$
generated by $\,\k\,$ (which is compact if $\,G\,$ is a finite
covering of a real form of a complex semisimple Lie group).
Since $\,\k$, $\,\p\,$ and  
$\,B\,$ are  $\,\Ad(K)$-invariant, 
 $\,\nu_m\,$ is right $\,K$-invariant, i.e. the
action of $\,G \times K\,$ on $\,G\,$ defined by 
$\,(g,\,k) \cdot l := glk^{-1}\,$
is by isometries. Here we allow discrete ineffectivity
given by the diagonal in
$\,Z(G) \times Z(G)$, where $\,Z(G) \subset K\,$ is the center of
$\,G$.  One has
$\,G \,= \,(G \times K)/H\,$ with 
$\,H\,$ the diagonal in $\,K \times K$.

Note that a different choice of a maximal compact connected
subalgebra $\,\k'\,$ induces an {\it equivalent} left-invariant Riemannian 
structure, i.e. there exists an isometric isomorphism
$$\,(G, \nu_m)  \to (G, \nu_m')\,.$$

\nmedskip
This is given by the internal conjugation transforming $\,\k\,$ in $\,\k'$. 

We summarize the main properties of the above
metrics in the following proposition where the 
degenerate metric $\,\nu_0\,$ can  be regarded  
as a limit case of non-degenerate ones.  

 
\nbigskip
\begin{prop}
\label{SEMISIMPLE} {\rm (\cite{HaIa2}, Sect.\,3, cf. \cite{Go}, proof of Thm.\,5.2)} $\ $Let $\,G\,$ be
a  non-compact, semisimple Lie group and, for $\,m \in \R$, let
$\,\nu_m\,$ be the above defined left-invariant
metric. Then

 \begin{enumerate}
\item[i)] the action of $\,G \times K \,$ by left and right 
multiplication is by isometries and

\medskip
\item[ii)] 
the direct sum 
$\,\h \oplus \m$, with $\,\h\,$ the isotropy Lie algebra and
$$\m := \{\,(\,-mX_\k+X_\p,\, -(1+m)\,X_\k\,) \in \g \times \k \ : \ 
X_\k+ X_\p \in \k \oplus\p \,\}\,,$$
is a naturally reductive decomposition of $\,\g \times \k\,=\, Lie(G\times K)$.
In particular for every 
$\,X= X_\k+X_\p\,$ in $\,\k \oplus \p \cong T_eG\,$ the unique geodesic through $\,e\,$
and tangent to $\,X\,$ is given by $\,\gamma_X: \R \to G$,
$$t \ \longrightarrow \ \exp_G \,t\,(\,-m X_\k+X_\p\,)\,
\exp_K \,t\,(1+m) X_\k\,.$$
\end{enumerate}

\end{prop}

\bbigskip
From \cite{Go}, section$\,$5, it follows that in the Riemannian cases $\,m>0\,$
the connected component of the isometry group is essentially given
 by $\,G\times K\,$ (here discrete ineffectivity is allowed)
 while in the pseudo-Riemannian symmetric case $\,m=-1\,$ it 
coincides with $\,G\times G$. Further information regarding the
Levi Civita connections (or their  limit when $\,m\,$ vanishes), the curvature tensor, the scalar and Ricci curvature
can be found in \cite{HaIa2}, section 3.


Let $\,M\,$ be a complete real-analytic Riemannian
manifold. Following the results of  Guillemin-Stenzel (\cite{GuSt})
and Lempert-Sz{\H o}ke
(\cite{LeSz}) one can introduce 
a complex structure on a subdomain of the tangent bundle
$\,TM\,$ which is {\it canonically} adapted to the given 
Riemannian structure. Recently R. Bielawski (\cite{Bi}) and
R. Sz{\H o}ke (\cite{Sz3}) have pointed out
that for the existence of such a  complex structure it is enough 
to have a real-analytic linear connection $\,\nabla$.
A real-analytic complex structure
on a domain $\,\Omega\,$ 
of $\,TM\,$ is adapted to $\,\nabla\,$ if all leaves of the
induced foliation are complex submanifolds with their
natural complex structure. That is,
for any $\,\nabla$-geodesic
$\,\gamma:I \to M \,$ the induced map $\,\gamma_*: TI\subset \C \to TM\,$
defined by $\,(x + iy) \mapsto y \, \gamma'(x)\,$ is holomorphic
on $\,(\gamma_*)^{-1}(\Omega)\,$
with respect to the adapted complex structure.
Here $y \, \gamma^{\prime}(x) \in T_{\gamma(x)} M$
is the scalar multiplication in the vector space $T_{\gamma(x)} M$.

The adapted complex structure exists and
it is unique on a sufficiently small neighborhood of $\,M$,
 which is identified with the zero section in its tangent
bundle $\,TM$.
If $\,\Omega\,$ is a domain of $\,TM\,$ containing
$\,M\,$ on which this structure is defined, then
we refer to it as an {\it adapted complexification}.

Associated to every non-degenerate metric $\,\nu_m\,$ of the above introduced one-parameter family one has the Levi Civita connection. For $\, X\, $ and $\,Y \,$ in $\,\g\,$ this is
given by the formula (cf. \cite{HaIa2})
$$
\nabla_{m\, X} Y = \frac{1}{2} \left( \lbrack X, Y \rbrack 
          + (1+m) \left(  \lbrack  X_{\mathfrak{k}},Y_{\mathfrak{p}}
                        \rbrack + \lbrack  Y_{\mathfrak{k}} , X_{\mathfrak{p}}
              \rbrack  \right) \right).$$

\smallskip
\noindent
Note that this uniquely defines a left-invariant,
real-analytic, linear connection also in the degenerate case $\,m=0$.
Therefore for all real $\,m\,$  one has an adapted complex structure
at least in a neighborhood of $\,M\,$ in  $\, TM$.
 In the next section we will see that  
there exists an adapted complexification which is maximal in the 
sense of containing any other adapted complexification.


\nbigskip
\begin{rem}
\label{AGUILAR}
$\,${\rm Let $\,M\,$ be a real-analytic, Riemannian
manifold with a free  action by isometries  of a compact 
Lie  group $\,K\,$ and endow $\,M/K\,$ with the unique
Riemannian metric such that the natural projection
$\, M \to M/K\,$ becomes a Riemannian
submersion. Then, as a consequence of results in 
 \cite{Ag}, if the adapted complex structure exists
 on all of $\,TM$, so it does on $\,T(M/K)$.
 In the pseudo-Riemannian context the analogous
 result does not hold in this generality.
 
 For instance, let $\,U\,$ be a compact, semisimple Lie group and 
 denote by $\,U^\C\,$ its universal complexification.
 Let  $\,G\,$ be  a non-compact, real
 form of $\,U^\C\,$ and $\,K = G\cap U$. Denote by $\,\k\,$ and $\,\u\,$ the Lie algebras of 
  $\,K\,$ and $\,U$, respectively, and by $\,B\,$ the Killing form on $\,U$.
 Consider the unique left-invariant Riemannian metric
 $\,\nu\,$ on $\,U\,$ defined for all $\,X,Y \in \u\,$ by
$$\nu|_{\u}(X,Y)=-2B(X_\k,Y_\k)- B(X_\p,Y_\p),$$ 
 
 \nmedskip
where  $\,\p:=\k^{\perp_B}$.
Endow $\,U \times K\,$
 with the unique bi-invariant pseudo-metric $\,\tilde \nu\,$
 such that
$$\tilde \nu|_{\u \times \k}((X,Z),(Y,W))=-B(X,Y)+ 2B(Z,W)$$
 
 \nsmallskip
 for all $\,(X,Z),\,(Y,W)\,$ in $\, \u\times \k$.
Then the projection 
 $\,U \times K\,\to\,U,\  (u,\,k) \to uk^{-1}\,$ turns out to be 
 a pseudo-Riemannian submersion. 
 Moreover the adapted complex
 structure is defined on all of $\,T(U\times K)$. Indeed $\,(U \times K, \tilde \nu) \,$
 is essentially (up to the sign of the metric in the second component) the product
 of two symmetric Riemannian spaces of the compact type, thus this
 is a consequence of results in \cite{Sz1} and \cite{Sz3}.
 
 However $\,(U,\nu)\,$ has some negative sectional
 curvatures (see \cite{DZ}, cf. \cite{HaIa2}, Sect.\,3),
 thus by Theorem 2.4 in \cite{LeSz} the adapted complex structure
 is not defined on all of $TU$.} 
 \qed
\end{rem} 
 
\medskip
\bigskip

\section{A family of maximal adapted complexifications}

\bigskip
Let $\,G\,$ be a  connected, non-compact, semisimple Lie group and
consider the one-parameter family of left-invariant 
metrics (pseudo-Riemannian for $\,m<0$, degenerate for $\,m=0$)
introduced in  section 2 and defined by 
$$\nu_m|_\g :=  -mB(X_\k,Y_\k) + B(X_\p,Y_\p).$$

\nsmallskip
Then (cf. Prop.$\,$\ref{SEMISIMPLE}) the action of $\,L= G \times K\,$
 by left and right multiplication is by isometries,
the isotropy in $\,e\,$ is $\,H= \{\,(k,k) \in G\times K \,:\, k \in K \}\,$
and the quotient $\,L/H\,$ is a natural reductive realization of $\,(G, \nu_m)$.
The Riemannian exponential map $\,\Exp_e\,$ in $\,e\,$ is given by 
$$\Exp_e(X) = \exp_L(-mX_\k + X_\p, \, -(1+m)X_\k) \cdot e=$$ 
$$\exp_G(-mX_\k + X_\p)\exp_K(1+m)X_\k \, ,$$

\medskip
\noindent
for every $\,X \in \g \cong T_eG$. Note that the $\,L$-action on $\,G\,$
induces an action on the tangent space $\,TG\,$ just by differentiation.
If one identifies $\,TG\,$  with $\,G \times \g\,$ as usual,
this action reads as
$$\,(g,k) \cdot(g', X)=(gg'k^{-1}, \, \Ad_k(X)).$$
The group $\,L\,$ also acts on the universal complexification $\,G^\C\,$ of $\,G\,$
by left and right multiplication.

Next we show the existence of a maximal
complexification in $\,TG\,$ adapted to the connection associated to 
$\,\nu_m\,$ for all real $\,m$. It can be characterized as follows.


\nbigskip
\begin{prop}
\label{POLARMAP} Let $\,G\,$ be a connected, non-compact, semisimple Lie group
endowed with a metric $\,\nu_m\,$ of the above family.
Then there exists a maximal  complexification
$\,\Omega_m\,$ adapted to the connection associated to 
$\,\nu_m$.
Let  $$P_m: G \times \g \to G^\C$$
be the $\,L$-equivariant map defined by 
$$ (g, X) \to g \exp_{G^\C}i (-mX_\k + X_\p) \exp_{K^\C}
i(1+m)X_\k \,.$$

\smallskip
\noindent
Then $\,\Omega_m\,$ is given by 
the connected component of $\,\{\,|DP_m| \not= 0\,\}\,$
containing $\,G \times \{0\}$. The polar map
$\,P_m|_{\Omega_m}\,$ is locally biholomorphic.
\end{prop}

\smallskip
\begin{proof}
Note that the universal complexifications of $\,L\,$ and $\,H\,$ are 
$\,L^\C=G^\C \times K^\C\,$ and
$\,H^\C = \{\,(k,k) \in L^\C \,:\, k \in K^\C \}$, respectively.
Let $\,L\,$ act on $\,L^\C/H^\C\,$ by left multiplication and
consider the $\,L$-equivariant map 
$\,\tilde P_m: TG \to L^\C/H^\C \,$ defined, for $\,l \in L\,$ and $\,X \in T_eG$, by 
$$l_*(X) \to l \exp_{L^\C}(i (-mX_\k + X_\p, \, -(1+m)X_\k))H^\C\,.$$

\nsmallskip
Identify $\,G=L/H\,$  with the zero section in $\,TG$.
Then an analogous argument as in \cite{HaIa1},$\,$Corollary$\,$3.3, applies to show that
the connected component $\,\Omega_m \,$ of $\,\{\,|D\tilde P_m| \not= 0\,\}\,$
which contains $\,G\,$ is  the maximal adapted complexification 
and the restriction $\,\tilde P_m|_{\Omega_m}\,$ 
is  locally biholomorphic.

Since $\, \exp_{L^\C} =\exp_{G^\C} \times \exp_{K^\C}$,  the statement is a consequence of the following real-analytic, $\,L$-equivariant identification
$$    L^\C/H^\C \to G^\C,  \quad   (g,k)H^\C \to gk^{-1}\ .$$
\end{proof}

 
\bigskip
In order to  describe $\,\Omega_m\,$
it  is convenient to determine a slice for the $\,L$-action.
Let $\, \a^+\,$ be the closure of  a Weyl chamber in a maximal 
abelian subalgebra $\,\a\,$ of $\,\p\,$ and define 
$$ \Sigma := \{\,(e,X)\in G  \times \g \ : \ X_\p \in \a^+ \}\,.$$

\nsmallskip
Since $\,\a^+\,$ is a fundamental domain for the $\,\Ad_K\,$ action on
$\,\p$,  every $\,L$-orbit of $\,G  \times \g\,$ meets 
$\,\Sigma$. Then one has


\bigskip
\begin{prop}
\label{SLICE}
Denote by $\,\Sigma_m\,$  the connected component  of  $\,(e,0)\,$
in the subset $\,\{\, (e,X) \in \Sigma \ : \   |(DP_m)_{(e,X)}| \not=0 \,\}\,$
of $\,\Sigma$. 
Then $\,\Omega_m = L \cdot \Sigma_m$.
\end{prop}

\smallskip
\begin{proof}
 Note that the $\,L$-equivariance of $\,P_m\,$ induces that
 of $\,DP_m$, i.e. 
$$(DP_m)_{l\cdot (e,X)} \circ Dl_{(e,X)}=  Dl_{P_m(e,X)} \circ (DP_m)_{(e,X)}\,,$$

\nsmallskip
for all $\,l \in L\,$ and $\,(e,X) \in \Sigma$.
In particular $\,(DP_m)_{l\cdot (e,X)}\,$ has maximal rank if and only if so does  
$\,(DP_m)_{(e,X)}$, implying the statement. 
\end{proof}

\bigskip

\section{The case of $\,SL_2(\R)$}

\bigskip
Let $\,G=SL_2(\R)$.
Here we choose $\,K=SO_2(\R)\,$ and give an explicit description of $\,\Sigma_m\,$
in terms of  fixed basis of $\,\k,\, \a,\, \p$. By Proposition
\ref{SLICE} this determines the maximal adapted complexification 
$\,\Omega_m\,$ associated to the connection of $\,\nu_m$.
Identify  $\,{\mathfrak sl}_2(\R)\,$ with the set of  zero trace
matrices
and let 
$$U =\left ( \begin{matrix}  \,0  &  -1\,   \cr
                                     \, 1  & \ 0   \,\cr 
\end{matrix} \right ), \quad 
H =\left ( \begin{matrix}  \, 1  &  \ 0\,   \cr
                              \, 0  &    -1 \,  \cr 
\end{matrix} \right ), \quad 
W =\left ( \begin{matrix} \, 0 \  & \   1 \,  \cr
                                \,1 \ &  \  0  \,  \cr 
\end{matrix} \right ).$$

\nsmallskip
Then $\,\{\,U\,\}\,$ is a basis of $\,\k$, while $\,\{H, W\}\,$ gives a basis of $\,\p$.
Consider the polar map $\,P_m: G \times \g \to G^\C\,$ introduced
in section 3 and choose
$\,\Sigma\,$ by fixing $\, \a^+:=\{\,aH \ : a \ge 0 \,\}$.
Given $\,(e,X) \in \Sigma\,$ consider the natural identification
$\,T_{(e,X)}(G \times \g) \cong \g \times \g\,$
and note that for $\,Y \in \g\,$ one has   

\begin{align*}
(DP_m&)_{(e,X)}(Y,0)=\frac{d}{ds}\bigr|_0 P_m(\exp sY,X) =
\frac{d}{ds}\bigr|_0 \exp sY \exp i (-mX_\k + X_\p)  \\
& \quad \exp i(1+m)X_\k=DR_{\exp i(1+m)X_\k} \circ DR_{\exp i (-mX_\k + X_\p)} Y,\\
& \\
(DP_m&)_{(e,X)}(0,Y_\k)=\frac{d}{ds}
\bigr|_0 P_m(0,X +sY_\k) =\frac{d}{ds}\bigr|_0 \exp i (-m(X_\k+sY_\k) + X_\p)\\
& \exp
i(1+m)(X_\k+sY_\k) = DR_{\exp i(1+m)X_\k} \circ (D\exp)_{i (-mX_\k + X_\p)} (-imY_\k) \ +\\
&  \quad \quad DL_{\exp i (-mX_\k + X_\p)} \circ (D\exp)_{i(1+m)X_\k}(i(1+m)Y_\k) ,\\
& \\
(DP_m&)_{(e,X)}(0,Y_\p)=\frac{d}{ds}
\bigr|_0 P_m(0,X +sY_\p) =\frac{d}{ds}\bigr|_0 \exp i (-mX_\k + X_\p+sY_\p) \\
&\quad \exp i(1+m)X_\k =DR_{\exp i(1+m)X_\k} \circ (D\exp)_{i(-mX_\k + X_\p)}(iY_\p) \,,\\
\end{align*} 

\noindent
where $\,L_g\,$ and $\,R_g\, $ denote left and
right multiplication by $\,g \in G^\C$, respectively, and $\,\exp\,$
is the exponential map of $\,G^\C$. Recall that by identifying 
$\,T_{\exp X}G^\C \,$ with $\,\g^\C\,$ via left multiplication one has
$$\,(D\exp)_X(Y) = \sum_{l=0}^\infty \frac{(-1)^l}{(l+1)!}\ad ^l(X)Y$$

\smallskip
\noindent 
for all $\,X,Y \in \g^\C\,$ (see, e.g. \cite{Va}). 
Fix $\,X=uU +aH \,$ in $\, \Sigma$. By using the above formulae, one shows that  the 
basis $\,\{(U,0), \,(H,0),\,(W,0),\,(0,U), \,(0,H),\,(0,W)\,\}$ of
$\, \g \times \g  \cong T_{(e,X)}(G\times \g)\,$
is mapped by the differential $\,DP_m|_{(e,X)}\,$ into the image via  $\,\Ad_{\exp- i(1+m)X_\k}\,$ of the following six vectors

$$\Ad_{\exp i(umU-aH)}U, \ \ \Ad_{\exp i(umU-aH)}H, \ \ 
\Ad_{\exp i(umU-aH)}W,$$

$$ \sum_{l=0}^\infty \frac{(-1)^l}{(l+1)!}\ad^l(i(-umU+aH))(-imU) \ +  
i(1+m)U,$$

$$\sum_{l=0}^\infty \frac{(-1)^l}{(l+1)!}\ad^l(i(-umU+aH))(iH), \ \  
\sum_{l=0}^\infty \frac{(-1)^l}{(l+1)!}\ad^l(i(-umU+aH))(iW).$$

\nmedskip
Now note that  $\,w \to \cosh w\,$ and $\,w \to \frac{\sinh w}{w}\,$ are even holomorphic,  thus the functions $\,z \to \cosh \sqrt z\,$  and $\,z \to \frac{\sinh \sqrt{z}}{\sqrt z}$
are well defined and holomorphic. Moreover the cofficients of their Taylor
series around the origin are real.
Then, by restriction one obtains  two real-analytic functions
which will be denoted in the sequel by $\,C\,$ and $\,S$, respectively.




\nsmallskip
 Letting
$\,x=4u^2m^2 -4a^2$, 
a further computation (see Appendix) shows that the above six vectors
can be written as

\smallskip
$$\left ( 1-4a^2\frac{C(x)-1}{x} \right) U \ + \ 
 4aum\frac{C(x)-1}{x} \, H\  +\ 
2aS(x) \, iW\,,$$

\smallskip
$$-4aum\frac{C(x)-1}{x} \, U \ +\  
\left ( 1+4u^2m^2\frac{C(x)-1}{x} \right) H\  + \
2umS(x) \,iW\,,$$

\smallskip
$$2aS(x) \,iU \ -\
2umS(x) \,iH \ +\  
C(x)\, W\,,$$

\smallskip
$$\left ( 1+4a^2m\frac{S(x)-1}{x} \right) \,iU\  - \ 4aum^2\frac{S(x)-1}{x}\,iH \ + \ 
2am\frac{C(x)-1}{x} \,W\,,$$

\smallskip
$$-4aum
\frac{S(x)-1}{x}\, iU \ +  \ 
\left ( 1+4u^2m^2\frac{S(x)-1}{x} \right) \, \,iH \ -\  2um\frac{C(x)-1}{x} \,W\,,$$

\smallskip
$$-2a\frac{C(x)-1}{x} \,U \ + \ 
2um\frac{C(x)-1}{x}\,H\ +\  S(x)\,iW\,,$$

\bigskip
Let us point out without proof some properties of
the functions
 $\,C:\R \to \R$, $\,x \to \cosh \sqrt x $,
and $\,S:\R \to \R\,$, $\,x \to \sinh \sqrt x /\sqrt x$,
which are used in the sequel.


\bigskip
\begin{lem}
\label{FUNCTIONS}
$\,$For $\,x>-\pi^2\,$ the real-analytic functions $\,S, \,S'\,$ are
strictly positive and $\,S,\,C,\,S',\,C/S,\,S/S'\,$ are strictly
increasing. Moreover
$$C'(x) = \frac{1}{2}S(x)\,, \ \  \quad S'(x)=\frac{C(x)-S(x)}{2x}  \ \ \quad
 and  \quad \ \ x< \frac{C(x)}{2S'(x)}.$$
\end{lem}

\bigskip
We can now determine the maximal adapted complexification 
$\,\Omega_m\,$ by computing the slice $\,\Sigma_m\,$
introduced in Proposition \ref{SLICE}. 


\nbigskip
\begin{prop}
\label{DOMAIN} 
The slice  $\,\Sigma_m\,$ for the maximal adapted complexification
 consists of the elements $\,(e,uU+aH)\,$ in $\,\Sigma\,$
such that

\begin{itemize}
\smallskip
\item[$*)$]  
$\ \ 4u^2m^2 - 4a^2 > -\pi^2$
 
\medskip
\item[$**)$] $\ \  4u^2m^2 +m4a^2  < f(4u^2m^2 - 4a^2)\,,$
\end{itemize} 

\nsmallskip
where the function $\,f:\R \to \R\,$
is defined by $\,f(x) := \frac{x\, C(x)}{C(x)-S(x)}= \frac{C(x)}{2S'(x)}\, $.
\end{prop}

\smallskip
\begin{proof}
Let $\,X=uU+aH \in \Sigma\,$ and $\,x=4u^2m^2 - 4a^2\,$.  By the above computations
$\,(DP_m)_{(e,\,X)}\,$
is non-singular if and only if the two determinants  

$$\left | \begin{matrix}
1-4a^2\frac{C(x)-1}{x} &
 4aum\frac{C(x)-1}{x} &
2aS(x) \cr
\cr
-4aum\frac{C(x)-1}{x} &
1+4u^2m^2\frac{C(x)-1}{x} &
2umS(x) \cr
\cr
-2a\frac{C(x)-1}{x}   &
2um\frac{C(x)-1}{x} &  S(x) \cr
 \end{matrix} \right | $$

\noindent 
and 
 
$$\left |  \begin{matrix}
2aS(x)  & -  2umS(x)  &  C(x)  \cr
\cr
1+4a^2m\frac{S(x)-1}{x}  & 
-  4aum^2\frac{S(x)-1}{x} & 
2am\frac{C(x)-1}{x} \cr
\cr
-4aum\frac{S(x)-1}{x} &
1+4u^2m^2\frac{S(x)-1}{x} &
 -  2um\frac{C(x)-1}{x} \cr
 \end{matrix} \right |\,,$$

\nbigskip
do not vanish. A straightforward computation yields
the two conditions above.
\end{proof} 

\bigskip

Note
that in the pseudo-Riemannian symmetric case
$\,m=-1$, conditions $\,*)\,$ and $\,**)\,$ coincide
and yield the set $\,\{\,(e,X) \in \Sigma \, : \, D\exp_X \, {\rm not} 
\, {\rm singular}\,\}$. 


\nbigskip
\begin{defi}$\ $
For $\, m \in \R\,$ let $\,\Sigma_m^*\,$ be the subdomain of $\,\Sigma\,$ defined by 
condition $\,*)\,$ in Proposition \ref{DOMAIN}, i.e.

$$\Sigma_m^*:= \{ \,(e,uU+aH) \in \Sigma \ :\   4m^2u^2 -4a^2> -\pi^2 \,\}.$$
\end{defi}


\nmedskip
\begin{rem} 
\label{THESAME} 
{\rm (see Picture 1) $\ $For $\, m \leq -1\,$ one has
$\, \Sigma_m = \Sigma_m^*$. Indeed 
by Lemma  \ref{FUNCTIONS} one has  $\,x< f(x)$,
for $\, x >-\pi^2$.
Consequently for any $\,(e,uU+aH)\in \Sigma^*_m\,$
\smallskip
 $$ 4u^2m^2 +4a^2m \leq 4u^2m^2 -4a^2 <
  f(4u^2m^2 -4a^2),$$
  
\noindent 
\smallskip
thus condition $\,**)\,$ is automatically fulfilled. 

For $\,m > -1\,$
the closure  of   $\,\Sigma_m\,$
in $\,\Sigma\,$ is contained in 
$\,\Sigma_m^*$. Indeed, for 
$\, 4u^2m^2 -4a^2 = -\pi^2$, one has  
$$ 4u^2m^2 +4a^2m >4u^2m^2 -4a^2 = -\pi^2 =f(-\pi^2)
=f(4u^2m^2 -4a^2)\,,$$

\nsmallskip
therefore condition $\,**)\,$ is not fulfilled.
In particular the boundary of
$\, \Sigma_m\,$ in $\,\Sigma\,$  is defined by
$$ \{ \,(e,uU+aH)\in \Sigma^*_m \ :\  
4u^2m^2 +4a^2m=f( 4u^2m^2 -4a^2) \,\}.$$}
 \qed 
\end{rem}


\begin{figure}[h]
    \leavevmode 
    \centerline{Picture 1}
  \centerline{$\ $}
    \centerline{$m \,\le \, -1 \quad\quad\quad\quad\quad\quad \quad 
           \ \      \ \    m\, =\, -\frac{1}{2}  \ \ \ \ 
\quad\quad\quad\quad\quad\quad \quad \quad m\,=\,1$}
 \begin{center}
    \epsfxsize=0.34\textwidth     
    \epsffile{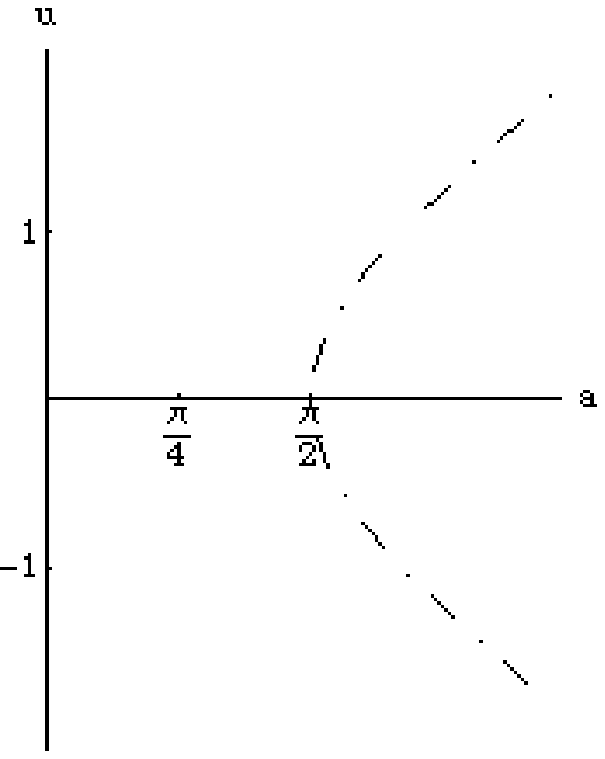}\ \ \  \quad
    \epsfxsize=0.225
    \textwidth\epsffile{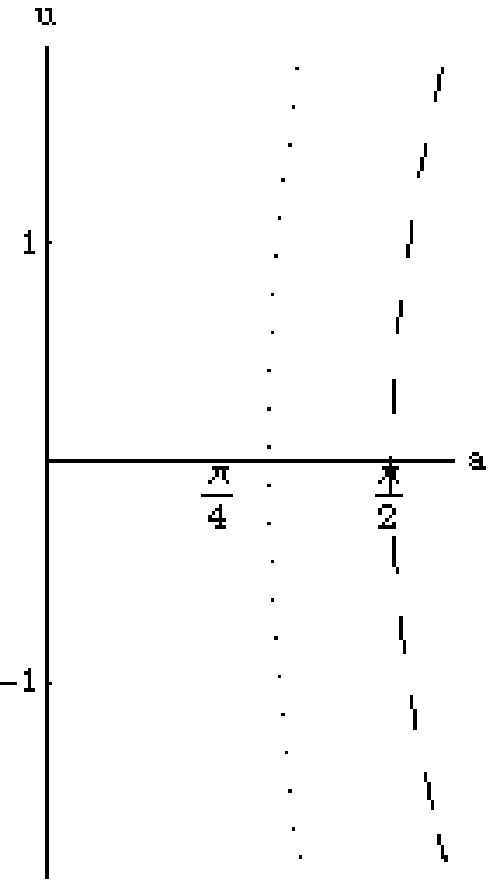}  \ \ \ \quad
     \epsfxsize=0.3
     \textwidth
     \epsffile{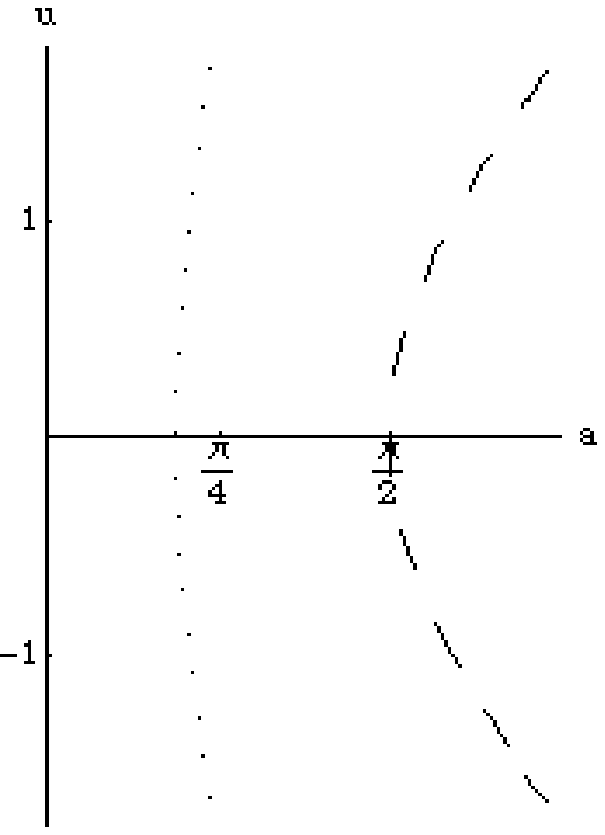}
      \end{center}
  \centerline{ $\cdot \, \cdot \, \cdot \,$ boundary of $\,\Sigma_m
\quad \quad \quad \quad \quad \quad$ - - - 
  $\,$boundary of $\Sigma_m^*$}
\end{figure} 

\nmedskip
\begin{rem}$\ $ 
{\rm For $\,m>-1\,$ one checks that the vector fields tangential to
$\, P_m(\Sigma)$, i.e. $\,(DP_m)_{(e,X)}(0,U)\,$ and
 $\,(DP_m)_{(e,X)}(0,H)$, remain 
linearly independent on the boundary of $\,P(\Sigma_m)$. 
Thus the boundary of the maximal adapted complexification
can be characterized by saying that 
$\,L$-orbits (in fact $\,\exp(W)$-orbits)
become tangential to 
$\,P_m(\Sigma)$.

On the other hand for $\,m\leq -1\,$ 
it is the dimension of  $\,L$-orbits at the boundary of
$\,P(\Sigma_m)\,$ to drop from 4 to 3 
(see $\,i)\,$ of Prop. 
\ref{PROPERTIES} below), giving 
again a characterization of the maximal
adapted complexification.}
 \qed
\end{rem}

\bigskip
Before studying the restriction of $\,P_m\,$ to $\,\Sigma_m\,$ it will be
useful to obtain a concrete realization of quotients of $\,G^\C\,$ with 
respect to those actions which are involved.

\medskip
\bigskip

\section{Slices and quotients of $\,SL_2(\C)$}

\bigskip
Here  $\,G\,$ and $\,K\,$ denote
$\,SL_2(\R)\,$ and $\,SO_2(\R)$, respectively. Let
$\,L=G\times K\,$ act on $\,G^\C\,$ by left and right multiplication.
Note that since $\,K\,$ is compact, this action is proper.
The main goal of this section is to present models for
the quotients $\,G \backslash G^\C$, $\,G^\C/L,$ $\,G^\C/(G \times K^\C)\,$ 
and for the relative quotient maps. 
First consider the map
$$\,\Pi_1:G^\C \to G^\C\, , \quad \,g \to \sigma_G(g)^{-1}g\,,$$

\smallskip
\noindent
 where $\,\sigma_G:G^\C \to G^\C$, $\, g \to \overline g$,  is the 
conjugation in $\,G^\C$, i.e. the
unique antiholomorphic involutive automorphism of $\,G^\C\,$
whose fixed point set is $\,G$.
Let $\,G\,$ act on $\,G^\C\,$ by left multiplication and
note that every fiber of this map consist of a single $\,G$-orbit.
 Thus $\,\Pi_1(G^\C)\,$ is set theoretically 
equivalent to $\, G\backslash G^\C\,$ and
$$\Pi_1:G^\C \to \Pi_1(G^\C)$$

\nsmallskip
is a realization of the quotient map.
Moreover, a simple computation shows that (cf., e.g. \cite{Zh})
$$\,\Pi_1(G^\C)=\{\,g \in G^\C \ : \ \sigma_G(g)=g^{-1} \,\}.$$

\smallskip
\noindent
It is  convenient to
consider the automorphism $\,A:G^\C \to G^\C\,$ transforming 
$\,SL(2,\R)\,$ onto $\,SU(1,1)$. This is induced by 
the unique complex Lie algebra morphism of $\,\g^\C\,$ mapping the 
basis $\,\{U,\,H,\, W\}\,$ (cf. beginning of 
Sect. 4) into $\,\{iH,\,iU,\, W\}$. Recall that
the involution of $\,G^\C\,$ defining $\,SU(1,1)\,$  is given by
 $\,\sigma_{SU(1,1)}(g) = J^t\overline g^{-1}J$,
where

$$\,J:=\left ( \begin{matrix}  1  &  \ 0   \cr
                                          0    & -1	  \cr 
\end{matrix} \right ).$$ 

\nsmallskip
Since the elements of the above basis are fixed by the Lie algebra
automorphisms induced by $\,\sigma_G\,$ and $\,\sigma_{SU(1,1)}\,$
respectively, it follows that 
\smallskip
\begin{equation}
\label{EQUI} 
A\circ \sigma_G = \sigma_{SU(1,1)} \circ A\,,
\end{equation}


\noindent
Note that $\, {\mathcal Q} := A \circ \Pi_1(G^\C)\,$
can be identified with $\,\Pi_1(G^\C) \cong G \backslash G^\C$.
Then 
\smallskip
$${\mathcal Q} = \{\, A(g) \ : \ g \in G^\C \  {\rm and} \ \sigma_G(g)=g^{-1}\,\}=
\{ \, g \in G^\C \, :\, \sigma_{SU(1,1)}(g)=g^{-1} \}=$$

$$\left \{  \left ( \begin{matrix}  \ s  &  b   \cr
                   -\overline b    &  t	  \cr \end{matrix} \right ) \, :
				   \,  s,t \in \R, \ b \in \C \, \ {\rm and}\, \
				   st +|b|^2=1 \right \}\,$$

\medskip
\noindent
gives a model of the quotient $\,G \backslash G^\C$.
Let us describe how the right $\,K^\C$-action on $\,G^\C\,$
is transformed after applying  $\,A \circ \Pi_1$. 
For $\,\lambda \in \C\,$ and $\,g \in G^\C\,$ one has 

$$A\circ \Pi_1(g \exp (-\lambda U)) = A(\sigma_G(g \exp(-\lambda U))^{-1}g\exp (-\lambda U))=$$

$$A(\sigma_G(\exp(-\lambda U)))^{-1}A(\sigma_G(g))^{-1}A(g)A(\exp (-\lambda U)))$$

\medskip
\noindent
and by (\ref{EQUI}) this gives 

$$\sigma_{SU(1,1)}(A(\exp(-\lambda U)))^{-1}A(\sigma_G(g)^{-1}g)A(\exp (-\lambda U)))=$$

$$(J\overline {\exp(-i\lambda H)}J)^{-1}A(\Pi_1(g)) \exp (-i\lambda H)=
\exp(i \overline \lambda H)A(\Pi_1(g)) \exp (-i\lambda H).$$

\medskip
\noindent
Thus $\,A\circ \Pi_1:G^\C \to {\mathcal Q}\,$ is $\,K^\C$-equivariant, if one let
$\,K^\C\,$ act on $\,G^\C\,$ by right multiplication and on
$\,{\mathcal Q}\,$ by 

$$\exp (\lambda U) \cdot \left ( \begin{matrix}  \ s  &  b   \cr
  -\overline b    &  t	  \cr \end{matrix} \right ) :=
  \exp (i\overline \lambda H)\left ( \begin{matrix}  \ s  &  b   \cr   
  -\overline b    &  t	  \cr \end{matrix} \right ) \exp (-i\lambda H)=$$

	$$\left ( \begin{matrix}  \ e^{2y}s  &  &  e^{2ix}b   \cr
\cr	
  -e^{-2ix}\overline b    &  & e^{-2y}t	  \cr \end{matrix} \right ),$$

\medskip
\noindent 
for every $\,x+iy = \lambda \in \C$.
In particular, after applying $\,A\circ \Pi_1$,
the right $\,K$-action on $\, G^\C\,$ reads as
rotations on $\,b$. Let 
\smallskip
$$\,{\mathcal P}:=\{\, (s,t) \in \R^2 \ : \ st\leq 1\}$$ 

\nsmallskip
and define 
$\,\Pi_2:{\mathcal Q} \to {\mathcal P}\,$ by
\smallskip
$$\left ( \begin{matrix}  \ s  &  b   \cr   
  -\overline b    &  t	  \cr \end{matrix} \right ) \to (s,t) \,.$$

\smallskip
\noindent
For every $\,(s,t) \in {\mathcal P}\,$ the inverse image 
$\,\Pi_2^{-1}(s,\, t)\,$ consists of 
a single $\,K$-orbit given by 
$$\left \{ \left ( \begin{matrix}  \ s  &  b   \cr   
  -\overline b    &  t	  \cr \end{matrix} \right ) \in {\mathcal Q} \ : \ |b|^2=1-st \right \}$$
  
\nsmallskip
In fact $\,{\mathcal P}\,$
is a realization of 
the quotient $\,{\mathcal Q}/K \cong G^\C/L$. Recall that 
 $\,L^\C= G^\C \times K^\C\,$ act 
 on $\,G^\C\,$ by left and right multiplication.
Let the one-parameter subgroup of $\,L^\C\,$
 defined by $\,R := \{e\} \times \exp i\k\,$ act on 
 $\,{\mathcal P}\,$ by $\, (\,\{e\} \times \exp(iyU)\,) \cdot (s,\,t) :=
 (e^{2y}s,e^{-2y}t)$, for all $\, y\in \R\,$
 and $\,(s,t) \in {\mathcal P}$. One has 
 

 \bigskip
\begin{lem}
\label{QUOTIENT1} Let $\,K^\C\,$ act on
$\,{\mathcal Q}\,$ by 

$$\exp (x+iy) U \cdot \left ( \begin{matrix}  \ s  &  b   \cr
  -\overline b    &  t	  \cr \end{matrix} \right ) =
  \left ( \begin{matrix}  \ e^{2y}s  &   e^{2ix}b   \cr
  -e^{-2ix}\overline b    & e^{-2y}t	  \cr \end{matrix} \right ),$$

\nsmallskip
for $\,x+iy \in \C$. Then 

\begin{itemize}
\smallskip
\item[$i)$]  a model for the quotient $\,G \backslash G^\C\,$ is $\,{\mathcal Q}\,$
with $\,K^\C$-equivariant quotient map $\,A\circ \Pi_1:G^\C \to {\mathcal Q}\,$. 
The  fixed point set  of the
$\,K$-action on $\,{\mathcal Q}\,$ is  given by $\,\{\left ( \begin{matrix}  s  &  0   \cr   
  0    &  t	  \cr \end{matrix} \right ) \ : \ st=1\}\,$,

\smallskip
\item[$ii)$] 
 a model for $\,{\mathcal Q}/K\,$ is given by $\,{\mathcal P}\,$
with quotient map $\,\Pi_2:{\mathcal Q} \to {\mathcal P}\,$
defined  by
$\,\left ( \begin{matrix}  \ s  &  b   \cr   
  -\overline b    &  t	  \cr \end{matrix} \right ) \to (s,t)\,$,

\smallskip
\item[$iii)$] 
a model for the quotient $\, G^\C/L\,$ is $\,{\mathcal P}\,$
with $\,R$-equivariant quotient map $\,F:=\Pi_2 \circ A \circ \Pi_1:G^\C \to {\mathcal P}$.
  In particular  $\,G \backslash G^\C/K^\C \cong {\mathcal P}/R$.
\end{itemize}
\end{lem} 


\medskip
\begin{rem}$\ $
\label{SIMMETRIES}
{\rm Let $\,N_{G^\C}(K^\C)\,$ denote the normalizer
of $\,K^\C\,$ in $\,G^\C$.
Then $\,e\,$ and $\,\tilde e:=\left ( \begin{matrix}  i  &  0   \cr   
  0    &  -i	  \cr \end{matrix} \right )\, $ represent the two elements
of the quotient $\,N_{G^\C}(K^\C)/K^\C$.
A simple computation shows that
$\,F\,$ is equivariant with respect to right (or left) multiplication by 
$\,\tilde e\,$ in $\,G^\C\,$ and reflection
 with respect to the origin in $\,{\mathcal P}$.
Furthermore $\,F\,$ is equivariant with respect to
conjugation by 
$\,\tilde e\,$ in $\,G^\C\,$ and reflection
 with respect to  the line of equation $\,\,{\,s=t\,}\,$ in $\,{\mathcal P}$ .}
 \qed
\end{rem}


\begin{figure}[b]
    \leavevmode 
    \centerline{Picture 2}
 \begin{center}
    \epsfxsize=0.47\textwidth     
    \epsffile{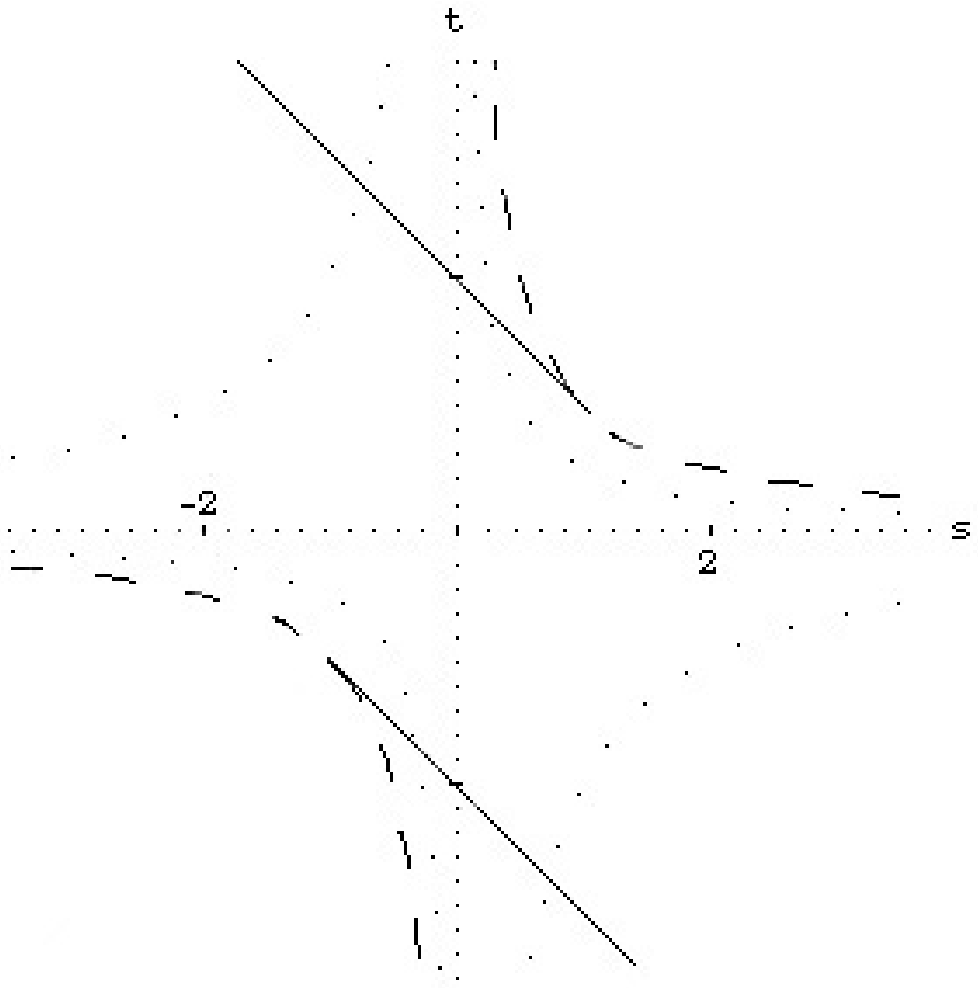}
      \end{center}
\medskip
   \centerline{
--- $\ \, {\rm  the\ slice\ {\mathcal S}}\quad \quad$
 - $\,$-  $\,$- ${\rm \ boundary\ of}\ {\mathcal P} \quad \quad$
$\cdot \, \cdot \, \cdot \ $ R-orbits}
 \centerline{$\ $}
\end{figure}



\noindent
\begin{rem}
\label{SLICEFORP}$\ $
{\rm The action of $\,R\,$ on $\,{\mathcal P}\,$ is not proper.
Set theoretically $\,{\mathcal P}/R\,$ can be identified
with the slice for the 
$\,R$-action  on $\,{\mathcal P}\,$  defined by the union
$\,\{\,s=t,\ st \le1\,\} \cup \,\{\,s=-t\,\}\cup \,\{\,p_1,\,p_2,\,p_3,\,p_4\,\}$,
where $\,p_1:=(2,0)$, $\,p_2:=(0,2)$, $\,p_3:=(-2,0)\,$ and $\,p_4:=(0,-2)\,$
correspond to the non-closed $\,R$-orbits in $\,{\mathcal P}$. The closure of these
orbits
is obtained by adding the unique fixed point $\,p_0=(0,0)$. However
for our purposes it is more convenient to consider the following slice
(cf. picture 2)
\smallskip
$$\,{\mathcal S}:=\{\,s+t= 2,\, t \ge s \,\} \cup \,\{\,s+t=-2 ,\, s \ge t\,\}\cup \,\{\,p_0,\,p_1,\,p_3\,\}.$$}
\qed
\end{rem}

\bigskip
Here and in the sequel a slice is assumed to intersect
every orbit in a single point.  
Let $\, \Sigma=\{ \,(e, X) \in G\times \g \ : \ X_\p \in \a^+\,\} \,$  be the slice
in $\,TG \cong G \times \g\,$ introduced in section 3. Consider the subdomain
$$\,\Sigma_{AG}=\{ \,(e, uU+aH) \in \Sigma \ : \  4u^2 -4a^2 > -\frac{\,\pi^2}{4} \,\}$$

\noindent
and denote by $\,\overline{\Sigma}_{AG}\,$ its closure in $\,\Sigma$.
One has


\bigskip
\begin{prop}
\label{PROPERTIES} $\,$Let $\, G\times K^\C \subset L^\C\,$ act on 
$\, G^\C\,$ by left and right multiplication. 
\begin{itemize}
 
\smallskip
\item[$i)$] 
 A slice for the $\,L$-action on $\,G^\C\,$
is given by 
\smallskip
$$S_1:=\exp i\overline{\Sigma}_{AG} \cup \tilde e \exp i\Sigma_{AG}.$$

\nsmallskip
All $\,L$-orbits are closed (the action is proper),
the union of all $\,3$-dimensional $\,L$-orbits is given by
$\,F^{-1}(\{\,st=1\,\})\,=L \cdot (\exp i\k\,\cup\,\tilde e \exp i\k)$.
All other orbits are $\,4$-dimensional with discrete isotropy given by
the ineffectivity $\,\pm(e,e)\,$ of $\,L$.

\smallskip
\item[$ii)$]  
A slice for the $\,G \times K^\C$-action
is given by 
$$\,S_2:=\{\,\exp \rho i(U+H) \,: \, \rho \ge 0\,\}\cup\tilde e
\{\,\exp \rho i(U+H) \, : \, \rho \ge 0\,\}\cup$$
$$
\{\,g_0,\,g_1,\,g_3\,\},$$

\nmedskip
where $\,g_0=\exp i \frac{\pi}{4}H$,
$\,g_1=\exp -\frac{i}{2}(U+H)\,$ and
$\,g_3=\tilde e g_1$.
The only non-closed $\,G \times K^\C$-orbits are those through
$\,g_1,\,g_2\,\,g_3,\,g_4$, with 
$\,g_2=\exp \frac{i}{2}(U+H)\,$ and
$\,g_4=\tilde e g_2$. Their closure is obtained by adding 
the orbit through 
$\,g_0$.
The only $\,4$-dimensional orbits are those through
$\,e,\,\tilde e$ and $\,g_0$, all other orbits have maximal dimension.
\end{itemize}
\end{prop}

\nmedskip
{\it Proof.}
For all $\,(e,uU+aH)\,$ in $\,\overline{\Sigma}_{AG}\,$ one has 
$$F(\,\exp i(uU+aH)) =  \Pi_2 \circ A (\exp 2i(uU+aH)) =
\Pi_2 (\exp -2(uH+aU)) =$$
$$\,\Pi_2 \left ( 
\begin{matrix} C(x) - 2uS(x)  &  2aS(x)  \cr
\cr
- 2aS(x)   &      C(x) + 2umS(x)  \cr
 \end{matrix} \right) =
 (C(x) - 2uS(x), \, C(x) + 2uS(x) ),$$
 
\nsmallskip
where $\,x=4u^2 - 4a^2$. Fix $\,x \ge -\pi^2/4\,$  and note that the set 
$\,Q_x:=\{\,(e,uU+aH) \in \overline{\Sigma}_{AG} \ : \ 4u^2 - 4a^2= x\,\}\,$
is given by $\,\{\,(e,uU+aH) \ : \  4u^2 \ge x, \ 0 \leq a= \sqrt{4u^2-x}/2\,\}.$
Then, from the above formula it follows  that  
$\,F(\exp iQ_x)\,$ consists of the intersection
of $\,{\mathcal P}\,$ with the line of equation $\,s+t=2C(x)$. As a consequence 
 $\,\exp i \overline{\Sigma}_{AG}$ is mapped  bijectively onto 
 $\,{\mathcal P} \cap \{s \ge -t\}$. This and  Remark \ref{SIMMETRIES} imply that $\,F\,$ maps $\,S_1\,$ bijectively onto $\,{\mathcal P}$. Thus $\,S_1\,$ is a slice for the $\,L$-action on $\,G^\C$.

Since $\,F=\Pi_2 \circ A \circ \Pi_1$, by $\,i)\,$ of Lemma \ref{QUOTIENT1} the only 
3-dimensional $\,L$-orbits are those through $\,F^{-1}(\{\,st=1\,\})\,$ and all others are 4-dimensional. Moreover
the isotropy of the $\,K$-action on $\,{\mathcal Q} \setminus \Pi_2^{-1}(\{\,st=1\,\})\,$
is given by the ineffectivity $\,\pm(e,e)$, implying that last claim in $\,i)$.

For $\,ii)\,$ note that $\,F(\,\exp i\rho(U+H))= (1-2\rho, 1+2\rho )\,$ (cf. the above
formula). This and Remark \ref{SIMMETRIES} apply to show that 
the restriction $\,F|_{S_2}:S_2 \to {\mathcal S}\,$ is bijective, where  
  $\,{\mathcal S}\,$ is the $\,R$-slice in $\,{\mathcal P}\,$ 
introduced in  Remark \ref{SLICEFORP}.
One has  a  commutative diagram of canonical quotients 
$$\begin{matrix} 
G^\C    & \buildover{ F}\longrightarrow &    G^\C/L  \cong {\mathcal P} \quad \quad 
\quad \quad  \ \  \cr
               &                    &  \cr
              & \searrow  & \downarrow \quad \quad \quad  \quad  
               \quad  \quad  \ \cr
                &                 &     \cr
              &    &   G^\C/(G \times K^\C)  \cong {\mathcal P}/R   \cr 
\end{matrix} $$

\nsmallskip
 where, by $\,iii)\,$ of Lemma \ref{QUOTIENT1}, the map $\,F\,$ is 
 $\,R$-equivariant. This gives a one to one correspondence
 between the  $\,G \times K^\C$-orbit space of $\,G^\C\,$ and 
 the $\,R$-orbit space of $\,{\mathcal P}$. As a consequence $\,S_2\,$ is a slice for the 
 $\,G \times K^\C$-action on $\,G^\C$. Moreover  the $\,G \times K^\C$-orbits
 through $\,g_0, \dots, g_4\,$ correspond to the five $\,R$-orbits
 in $\,\{st=0\}  \subset    {\mathcal P}$, implying the topological claim.
Finally, the dimension of every $\,G \times K^\C\,$ orbit can be obtained 
 by adding to the dimension of the corresponding $\,R$-orbit  the dimension of the $\,L$-orbit given in $\,i)$. Since every $\,R$-orbit different from the fixed point $\,(0,0)\,$ is 
 one-dimensional, this concludes the proof.
 \qed


\bigskip
Regard  $\,G/K\,$ as a Riemannian symmetric space of rank one.
Its maximal adapted complexification can be realized in
$\,G^\C/K^\C\,$ and it is usually described as
$\,\Omega_{AG}=G \exp (i[0,\frac{\pi}{4})H)K^\C$
(see \cite{AkGi}, \cite {BHH}).
 Its boundary is given by $\,\cup_{j=1}^3 Gg_jK^\C\,$ and we 
 refer to $\,Gg_0K^\C\,$ as its singular boundary. For later use 
 we note the following direct consequence of \cite{GeIa}, Theorem 6.1 and Example 6.3.
 

\bigskip
\begin{lem}
\label{ENVELOPE}
For  $\,G=SL_2(\R)$, let $\,\Omega\,$ be a $\,G$-invariant domain of $\,G^\C/K^\C\,$
which contains  $\,\Omega_{AG}\,$ and its
non-singular boundary, i.e. $\,\Omega\,$ contains the invariant subset $\,\overline \Omega_{AG} \setminus G \exp
(i\frac{\pi}{4}H) K^\C \,$. Then  $\, \Omega\,$ is holomorphically convex
if and only if it coincides with $\, G^\C/K^\C$.
\end{lem}

\bigskip

\section{The reduced polar map}

  \bigskip
  Let $\,G=SL_2(\R)\,$ be endowed with one of the metrics
  $\,\nu_m$.
  By Proposition \ref{SLICE} the associated maximal adapted complexification is given by
  $\,\Omega_m = L \cdot \Sigma_m$, where $\,\Sigma_m\,$ consists of the elements 
of  $\,\Omega_m\,$ which are in the slice $\,\Sigma\,$ for the $\,L$-action
  on $\,  G \times \g \cong TG$. Note that every $\,L$-orbit in  $\,G \times \g\,$
  intersects $\,\Sigma\,$ in a single element, thus $\,\Sigma \cong TG/L$. One has a commutative  diagram 
  
  $$\begin{matrix} \ \  TG \cong G \times \g     \     &  \buildover{P_m}  \to \   &  G^\C \cr
                  \cr    
              \downarrow\   &          &     \downarrow \,F\cr
\cr
       TG/L \cong \Sigma \ \       &    \buildover{\hat P_m}  \to \  & 
		\    G^\C/L \cong {\mathcal P}\,,\cr 
\end{matrix} $$

\nsmallskip
where $\,F:G^\C \to \mathcal P\,$ is the realization of the quotient map given  
in Lemma \ref{QUOTIENT1} and 
 $\,\hat P_m:=F \circ P_m|_\Sigma$.
Here we show that the polar map $\,P_m|_{\Omega_m}\,$ 
is injective if and only if so is the restriction
$\,\hat P_m|_{\Sigma_m} $. We also point out certain
properties of the map $\,\hat P_m\,$ which are 
used in the remaining sections in order to
discuss injectivity of $\,P_m|_{\Omega_m}$.
It will turn out that $\,\Omega_m\,$ is biholomorphic to an $\,L$-invariant domain 
of $\,G^\C\,$ if $\,m \leq 0$, while it is a non-holomorphically 
separable Riemann domain over $\,G^\C\,$ if $\,m>0$.
In both cases $\,\Omega_m\,$ is not holomorphically convex
and its envelope of holomorphy is biholomorphic to $\,G^\C$.


\bigskip
We first  compute the two components
of $\,\hat P_m\,$  with respect to the  basis $\,\{U,H\}$.
One has 

$$\hat P_m(e,uU+aH) = F \circ P_m(e,uU+aH)= F\bigl (\,\exp i(-umU+aH)
\exp iu(1+m)U\,\bigr )=$$

$$\Pi_2 \circ A \bigl (\, \exp iu(1+m)U \exp 2i \left (-umU+aH
\right ) \exp iu(1+m)U\, \bigr ) =$$

$$ \Pi_2 \bigl (\, \exp -u(1+m)H \exp 2  (umH-aU)
\exp -u(1+m)H\, \bigr ) =$$

\medskip
$$\,\Pi_2 \left ( 
\begin{matrix}
e^{-2u(1+m)}\left( C(x) + 2umS(x)\right)  &  2aS(x)  \cr
\cr
- 2aS(x)   &     e^{2u(1+m)}\left( C(x) - 2umS(x)\right)  \cr
 \end{matrix} \right),$$

\smallskip
\nbigskip
where $\,x=4u^2m^2 - 4a^2$. This gives 
 
 
 \begin{equation}
 \label{REDUCEDMAP}
\begin{matrix}
\hat P_m(e,uU+aH) = \cr
\cr
\left (\,e^{-2u(1+m)}\left( C(x) + 2umS(x)\right ),
\, e^{2u(1+m)}\left( C(x) - 2umS(x)\right)   \, \right).\cr
\end{matrix}
\end{equation}


\nmedskip
\begin{rem}
\label{RANK}$\ $
{\rm The map $\,\hat P_m\,$ has maximal rank on 
$\, \{ \,(e,\, X) \in \Sigma_m \ :\ X_{\a^+} \not= 0 \,\}$.
Indeed  it is easy to check that $\,F\,$ has maximal rank on
$\,F^{-1}(\{\,st\not=1\,\})$,
therefore so does $\,\hat P_m\,$ on the set $\,\Sigma_m \cap \hat P_m^{-1}(\{\,st\not=1\,\})$.
From  formula (\ref{REDUCEDMAP}) it follows that this set coincides with
$\,\Sigma_m \cap \{ \,(e,\, uU+aH) \ :\ a > 0 \,\}$.}
 \qed
\end{rem}

\bigskip
Recall that for all real $\,m\,$ the slice
$\,\Sigma_m\,$ is contained in the domain
$\,\Sigma_{m}^* = \{ \,(e,uU+aH) \in \Sigma \ : \   4m^2u^2 -4a^2> -\pi^2 \,\}\,$
(cf. Rem. \ref{THESAME}). For such bigger domain one has


\bigskip
\begin{lem}
\label{DIFFEO}
The restriction of $\,P_m\,$ to any  $\,L$-orbit
of $\, L  \cdot \Sigma_{m}^*\,$
is a diffeomorphism onto an $\,L$-orbit of $\,G^\C$.
\end{lem}

\nsmallskip
\begin{proof} 
Note that the isotropy of $\,L\,$ at a point $\,(e,X)\,$ of $\,\Sigma_{m}^*\,$
is given by $\,\{\,(k,k) \ : \ k \in K\,\}\,$ if $\,X \in \k$, it consists of the ineffectivity 
$\,\pm(e,e)\,$ otherwise. 
Moreover since $\,K\,$ and $\,\exp i\k\,$ commute,
the identity

$$(g,k) \cdot P_m(e,X) = P_m(e,X)$$

\nsmallskip
holds true if and only if 
$$g\exp i(-muU+aH)k =\exp i(-muU+aH)\,.$$

\nsmallskip
Now, a similar computation as in formula (\ref{REDUCEDMAP}) yields
$$\,F(\exp i(X))=
(\, C(x) + 2umS(x)\,
,\, C(x) - 2umS(x)\,)\,,$$
which belongs to 
$\,\{\,st=1\,\}\,$ if and only if $\,a=0$. Then,
as a consequence of  $\,i)\,$
of Proposition \ref{PROPERTIES}, the isotropy at 
$\,P_m(e, uU+aH)\,$ is given by
 $\,\{\,(k,k) \ : \ k \in K\,\}\,$ if $\,a=0\,$ or $\,\pm (e,e)\,$
otherwise, which proves the statement.
\end{proof} 


\bigskip
\begin{rem}$\ $ 
{\rm For $\,(e,\, uU+aH) \in \partial \Sigma^*_m\,$ one has 
$ \, 4u^2m^2- 4a^2=-\pi^2$.
Thus, by formula (\ref{REDUCEDMAP}) 
\nsmallskip
$$\,F \circ P_m(e, \,uU+aH)
=(-e^{-u(1+m)}, -e^{u(1+m)}) \in \{\,(s,\,t) \in {\mathcal P}\ : \ st=1\}\,.$$

\nbigskip
Then $\,i)\,$
of Proposition  \ref{PROPERTIES} implies that the dimension
of the $\,L$-orbit through
$\,P_m(e,uU+aH)\,$ is only three.
In particular an analogous statement as in the above lemma
does not hold  on domains larger than  $\,\Sigma_m^*$.}
 \qed
\end{rem}

\smallskip
\bigskip
Since $\, \Omega_m=L\cdot \Sigma_m\,$ and $\,\Sigma_m \,$ is contained in
$\,\Sigma_m^*$,
from the above lemma it follows that  $\,P_m|_{\Omega_m}\,$
is injective if and only if different $\,L$-orbits in $\,\Omega_m\,$
are mapped by $\,P_m\,$  to different $\,L$-orbits in $\,G^\C$.
Recalling that for $\,G=SL_2(\R)\,$ every orbit intersects $\,\Sigma_m\,$ in a single point and that 
$\,F:G^\C\to {\mathcal P}\,$ is a realization of the quotient map
with respect to the $\,L$-action on $\,G^\C$, one has

\bigskip
\begin{prop}
\label{INIINI}
The polar map $\,P_m|_{\Omega_m}\,$
is  injective if and only if $\,\hat P_m|_{\Sigma_m}\,$
is injective.
\end{prop}

\bigskip
We conclude this section with a technical result which will be repeatedly
used in the sequel.
Consider the two involutions
$$\, \alpha:\Sigma_m^* \to \Sigma_m^* \quad \quad X_\k + X_\p \to -X_\k + X_p$$
and  
$$\, \beta:{\mathcal P} \to {\mathcal P} \quad \quad (s,t) \to (t,s)$$

\nsmallskip
and denote by
$\,{\rm fix}(\alpha)=\{\,X\in \Sigma_m^* \ :
\ X_\k=0\,\}\,$ and $\,{\rm fix}(\beta)=\{\,(s,t)\in {\mathcal P} \ : s=t\, \}\,$
the associated fixed point sets.

Note that $\,\hat P_m|_{\Sigma_m^*}:\Sigma_m^* \to {\mathcal P}\,$
is equivariant with respect to these involutions.
As a consequence $\,\hat P_m({\rm fix}(\alpha))\,$ is contained 
in $\,{\rm fix}(\beta)$. Also consider the $\,\alpha$-invariant map 
$\,\Gamma:\Sigma_m^* \to \R\,$ defined by
$$\Gamma(uU+aH):=
\frac{m}{1+m}\,\frac {C(4u^2(1+m)^2)}{S(4u^2(1+m)^2)}-\frac{C(x)}{S(x)}\,,$$

\nsmallskip
with $\,x=4u^2m^2-4a^2$.

\bigskip
\begin{lem}
\label{FIXPOINTS}$\ $
Let $\, \alpha:\Sigma_m^* \to \Sigma_m^*\,$ and
$\, \beta:{\mathcal P} \to {\mathcal P}\,$ be the two
involutions defined above.
\begin{itemize}
\item[$i)$]  For $\,m \leq -1\,$ one has
$$\left (\hat P_m|_{\Sigma_m^*}\right)^{-1}({\rm fix}(\beta)) =
{\rm fix}(\alpha)\,.$$

\item[$ii)$]  For $\,m >-1\,$ one has {\rm(cf. picture 3)}

$$\left (\hat P_m|_{\Sigma_m^*}\right)^{-1}({\rm fix}(\beta)) =
{\rm fix}(\alpha) \cup {\rm graph}(\gamma)\,$$

\nmedskip
where $\,\gamma:\k  \to \a \,$ is the real-analytic map
implicitly defined by $\,\{\,\Gamma=0\,\}$.
\end{itemize}
\end{lem}

\begin{figure}[h]
    \leavevmode 
    \centerline{Picture 3}
    \centerline{$m\, =\, -\frac{1}{2}  \ \ \ \ 
\quad\quad\quad\quad\quad\quad
\quad\quad\quad \quad \quad  \quad \quad m\,=\,1
\quad \quad\quad \quad $}
 \begin{center}
    \epsfxsize=0.47\textwidth     
    \epsffile{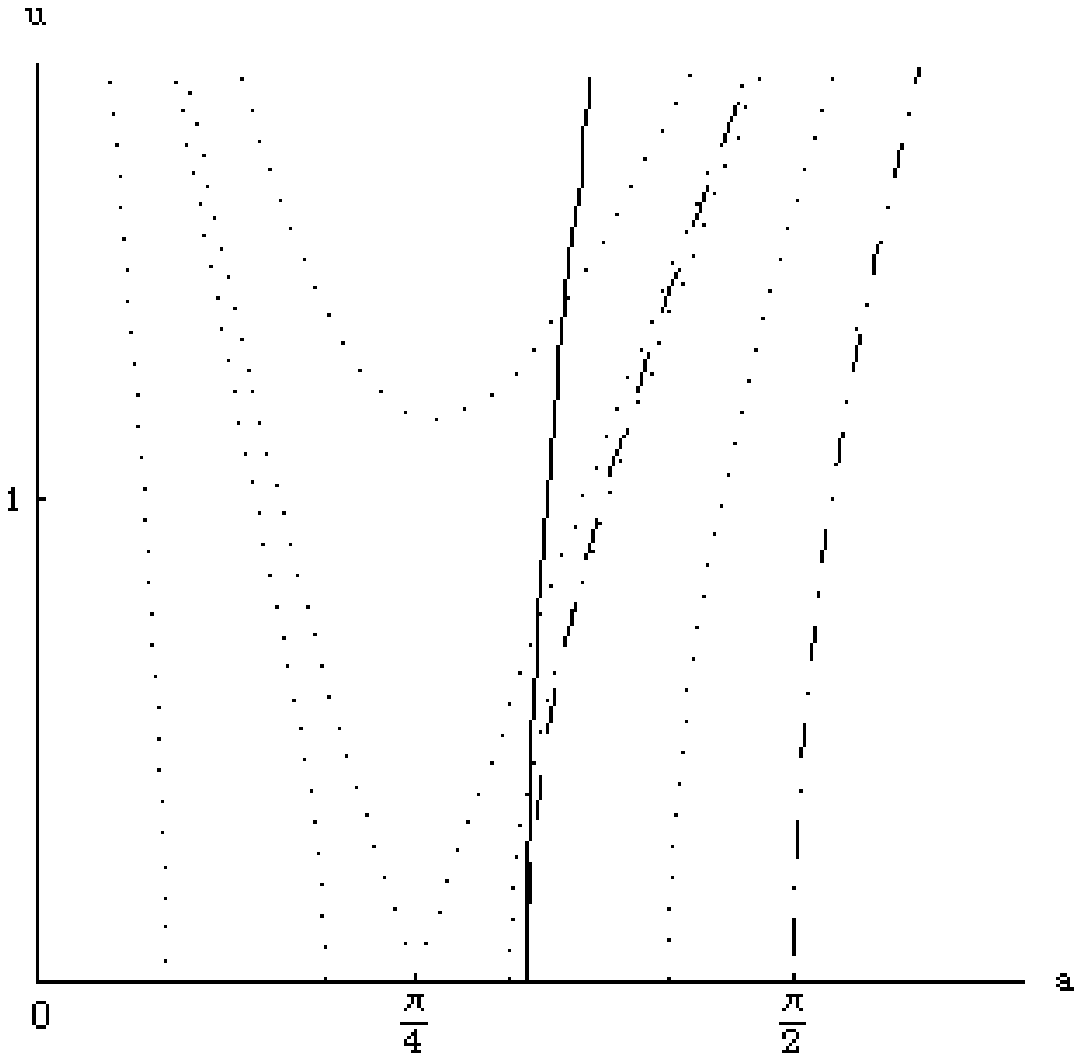} \quad 
    \epsfxsize=0.47
    \textwidth\epsffile{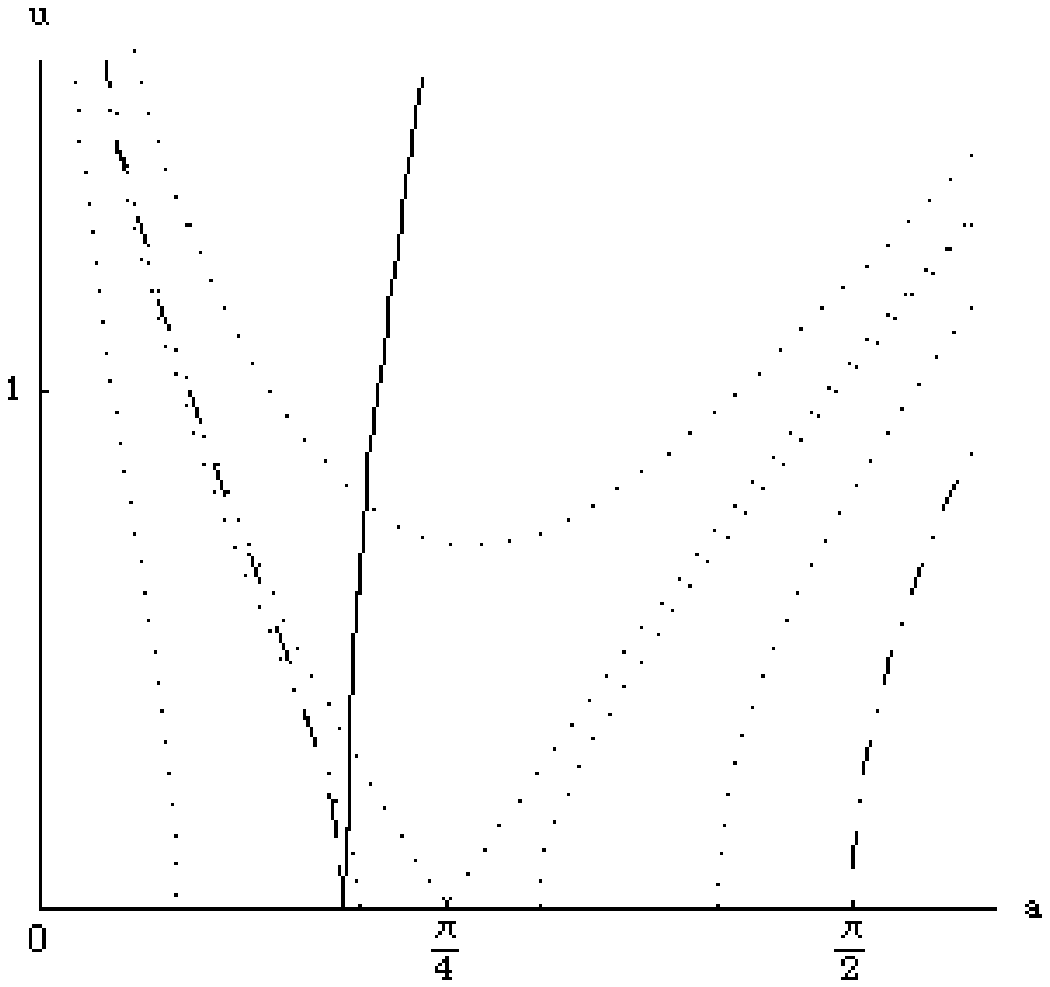} 
      \end{center}
\medskip
  \centerline{ --- $\ \, \partial \Sigma_m\quad $
- $ \cdot $ -  $\ \,\partial\Sigma_m^*\quad$
 - $\,$-  $\,$- $\ \{\,\Gamma =0 \,\}\quad$
$\cdot \, \cdot \, \cdot \ $ $\ell_c\ $ for 
$ \ c \approx 0.7,\ 0.2,\ 0,\ -0.7$}
 \centerline{$\ $}
\end{figure}

\bigskip
 \begin{proof}
Formula (\ref{REDUCEDMAP}) implies that the set 
 $\, \left (\hat P_m|_{\Sigma_m^*}\right)^{-1}({\rm fix}(\beta)) \,$ is given by
 
$$ \{\,uU+aH \in \Sigma_m^* \ : \    e^{-2u(1+m)}(\, C(x) + 2umS(x)\,)=
e^{2u(1+m)}(\, C(x) - 2umS(x)\,)   \,  \}$$

$$=\, \{\,  \cosh(2u(1+m))\, 2um \,S(x)  \,= \, \sinh(2u(1+m)) \,C(x)\, \}\, .$$

\medskip
\noindent 
Then the cases $\,m=-1\,$ and 
$\,m=0\,$ are straightforward. For $\,m\not=-1\,$ 
this set can be written as $\,\{\,u=0\,\}\cup\{\,\Gamma=0\,\}$, with 
$\,\Gamma\,$ as in the statement.

For $\,m<-1\,$ let  $\,u>0\,$ and note that 

$$\frac{\cosh (2u(1+m))}{\sinh(2u(1+m))}<
\frac{\cosh (2um)}{\sinh(2um)}.$$

\nmedskip
Recalling that
$\,x \to C(x)/S(x)\,$ is strictly increasing this yields  

$$\frac{m}{1+m}\frac {C(4u^2(1+m)^2)}{S(4u^2(1+m)^2)}
\,=\,\frac{\cosh (2u(1+m))\, 2um}{\sinh (2u(1+m))}\,>$$
$$\frac{\cosh (2um)\,2um}{\sinh (2um)}\,=\,
\frac {C(4u^2m^2)}{S(4u^2m^2)}\,\ge\,\frac{C(4u^2m^2-4a^2)}{S(4u^2m^2-4a^2)},$$

\smallskip
\noindent 
for any $\,a \geq 0$,
implying that $\,\{\,\Gamma=0\,\}$\, has no solution for $\,u>0$.
Along with the $\,\alpha$-invariance of   $\,\Gamma$, this implies $\,i)$.

For $\,m>-1$, $\,m\not= 0\,$ and $\,u>0\,$ fixed,
an analogous argument shows that 

$$\frac{m}{1+m}\frac {C(4u^2(1+m)^2)}{S(4u^2(1+m)^2)}<
\frac {C(4u^2m^2)}{S(4u^2m^2)}\,.$$

\smallskip
\noindent 
Since $\,C(x)/S(x)\,$ is strictly increasing for 
$\,x>-\pi^2\,$ and $\,C(x)/S(x) \to -\infty\,$
for $\,x\to  -\pi^2$, it follows that there exists a 
unique $\,a\in \R\,$ with $\,4u^2m^2-4a^2>-\pi^2\,$
such that $\,\Gamma(uU+aH)=0$. This and $\,\alpha$-invariance of   $\,\Gamma$
yield $\,ii)$.
 \end{proof}

\bigskip 

\section{The case of $\,m \leq -1$}

\bigskip
Let $\,G=SL_2(\R)\,$ be endowed with one of the metrics
  $\,\nu_m\,$ for some $\,m \leq -1$.
In this case we show that the polar map $\,P_m\,$ is a biholomorphism 
from the maximal adapted complexification $\,\Omega_m\,$ 
onto an $\,L$-invariant domain of $\,G^\C$. 
This domain is given by removing from $\,G^\C\,$ a
single $\,4$-dimensional $\,G \times K^\C$-orbit. Its 
envelope of holomorphy turns out to be biholomorphic to $\,G^\C$. 

\bigskip
 Consider the one parameter subgroup  $\,R:= \{e\}\times \exp(i\k) \,$ of 
 $\, L^\C$. The  $\{e\}\times K$-action
on $\,\Omega_m\,$ induces a local $\,R$-action  whose infinitesimal generator
 is given, for all $\,(g,\,X)\,$  in $\,\Omega_m$, by 
 $$\,iU \to J_m \bigl ( \, \frac{d}{dy}\big |_0(\{e\}\times \exp(yU))
\cdot (g,\,X) \,\bigr)\,.$$

\nsmallskip
 Here  $\, J_m\,$ denotes the adapted complex structure of $\,\Omega_m\,$
 
Since $\,P_m\,$ is holomorphic and $\,(\{e\} \times K)$-equivariant, it
is also locally equivariant with respect to such local 
$\,R$-action on $\,\Omega_m\,$ and the global $\,R$-action on $\,G^\C$.
Furthermore both (local) actions
commute with the $\,L$-actions on $\,\Omega_m\,$ and on $\,G^\C$,
thus they push down to (local) $\,R$-actions
on $\, \Omega_m/L  \cong \Sigma_m\,$
and on $\, G^\C/L \cong {\mathcal P}$, respectively.
Since one has the commutative diagram

$$\begin{matrix}  \Omega_m     \     &  \buildover{P_m}  \to \   &  G^\C \cr
                  \cr    
              \downarrow\   &          &     \downarrow \,F\cr
\cr
        \Omega_m/L \cong \Sigma_m       &    \buildover{\hat P_m}  \to \  & 
		{\mathcal P} \cong G^\C/L\,, \cr 
\end{matrix} $$

\medskip
\noindent
the restriction $\, \hat P_m|_{\Sigma_m}\,$ is locally $\,R$-equivariant.
Recall that by $\,iii)\,$ of Lemma \ref{QUOTIENT1},
the $\,R$-action on  $\,{\mathcal P}\,$ is explicitly given by 
$\,(\{e\} \times \exp iyU) \cdot (s,t) = (e^{2y}s,\, e^{-2y}t)$.
In particular the function $\,{\mathcal P} \to \R$, defined
by $\,(s,t) \to st$, is $\,R$-invariant. This is used to 
determine local $\,R$-orbits in $\,\Sigma_m\,$ as follows.
For  $\,c  \leq 1\,$ let $\,\ell_c:= 
(\hat P_m|_{\Sigma_m})^{-1}(\{\,st=c\,\})$
and denote by  $\,\Sigma_m^+\,$ the set
$\,\{\,(e, uU +aH) \in \Sigma_m \ : \ a>0, \ u>0\,\}$. 
One has
   

\bigskip
\begin{lem}
\label{LEVELSET}  Let $\,m\leq -1$. Then 
\begin{itemize}
\item[$i)$] a local $\,R$-orbit of $\,\Sigma_m^+\,$ coincides with a connected component
of $\,{\ell_c}\cap\Sigma_m^+$,

\nsmallskip
\item[$ii)$]  the intersection $\,{\ell_c}\cap\Sigma_m^+\,$ has two connected components
if $\,0 \le c < 1$, it is connected if $\,c<0$, 

\nsmallskip
\item[$iii)$]  $\,\hat P_m\,$ maps different local $\,R$-orbits of $\,\Sigma_m^+\,$ to
different $\,R$-orbits of $\,{\mathcal P}$,

\nsmallskip
\item[$iv)$]  $\,\hat P_m\,$ is injective on 
$\,\Sigma_m$,

\nsmallskip
\item[$v)$]  $\,\hat P_m\,$ maps different local $\,R$-orbits of $\,\Sigma_m\,$ to
different $\,R$-orbits of $\,{\mathcal P}$,

\nsmallskip
\item[$vi)$]  local $\,R$-orbits closed to $\,\{e\} \times \k\,$ are mapped 
bijectively by $\,\hat P_m\,$ to $\,R$-orbits of  $\,{\mathcal P}$, i.e.
$\,R\,$ acts globally in a $\,R$-invariant neighborhood of $\,\{e\} \times \k$.
\end{itemize}
\end{lem}

\nsmallskip
\begin{proof}$\ $
Note that $\,(0,0)\,$ is the unique fixed point for the $\,R$-action on
$\, {\mathcal P}\,$ and that $\,(\hat P_m)^{-1}(0,0)=(e, \frac{\pi}{4})\,$ does not belong
to  $\, \Sigma_m^+$.
Furthermore,  the restriction $\,\hat P_m|_{\Sigma_m^+}\,$
is locally diffeomorphic by Remark \ref{RANK} and locally $\,R$-equivariant by construction.
This implies that every local 
$\,R$-orbit of $\,\Sigma_m^+\,$ is one-dimensional.

One checks that $\, \ell_c = \{\,(e,\,uU +aH)
\in \Sigma_m^*\ :\ \phi_m(a,u)=\sqrt{1-c}  \, \}$,
where $\,\phi_m(a,u):=2aS(4m^2u^2-4a^2)$. Moreover
$ \partial \phi_m/\partial a \not= 0\,$ on $\,\Sigma_m^+$,
therefore $\, \ell_c\cap   \Sigma_m^+\,$ is a one-dimensional 
manifold. By construction, $\, \ell_c\,$ is locally $\,R$-invariant,
hence local $\,R$-orbits of $\,\Sigma_m^+\,$ are 
open and closed in $\, \ell_c\cap   \Sigma_m^+$, implying $\,i)$.

For $\,ii)$, recall that  $\,\Sigma_m=\Sigma_m^*\,$ by Remark \ref{THESAME}.
As a consequence $\,\phi_m\,$ vanishes on the boundary of
$\,\Sigma_m$. Also note that $\,\phi_m(a,0)= \sin (2a)\,$ and that
 for $\,u\,$ fixed and $\,a\,$ such that  $\, (e,\, uU +aU) \in \Sigma_m^+$,
the map 
$\, \R^{\ge 0} \to \R$, $\,u \to \phi_m(a,u)\,$ is strictly increasing
(cf. Lemma \ref{FUNCTIONS}).
This  implies that  $\, \ell_c\cap   \Sigma_m^+\,$
is the (connected) graph of a function defined on $\,\a^+ \setminus \{0\}\,$
for $\, c<0$, while, for  $\, 0 \le c < 1$, it consists
of two connected components 
(which are contained in $\,\Sigma_m^+ \setminus
\{\,(e,uU+aH)  \ : 0<a< \frac{\pi}{2}\ {\rm and} \ \sin (2a) \ge \sqrt{1-c} \ \}$).

For $\,iii)\,$ first consider the case
of $\, 0 <c < 1$. One needs to show that the two components of
$\, \ell_c\cap \Sigma_m^+\,$ are mapped to different components
of the hyperbola in $\, \mathcal P\,$ defined by $\,\{\,st=c\,\}$.
For this it is enough to note that this holds for the two limit 
points given by $\,\{\,(e,aH) \ : \ 0<a< \frac{\pi}{2}\ {\rm and} \ \sin (2a) =\sqrt{1-c} \
\}$.

The only other non-connected case is 
 $\, c =0$, when the  two
components of $\, \ell_0\cap   \Sigma_m^+\,$ have the same 
limit point $\,(e,\frac{\pi}{4}H).$ By Remark \ref{RANK}
the map $\,\hat P_m\,$ is  a diffeomorphism in a neighborhood (in $\,\Sigma_m$)
of this point, thus these two 
components are mapped to different components
of $\, \{\,(s,t) \in \mathcal P \ : \  st=0, \ (s,t)\not=(0,0)\,\}$, as claimed.
The case $\,c<0\,$ is straightforward.

Since all $\,R$-orbits in  $\,{\mathcal P}\,$
have connected isotropy and $\,\hat P_m|_{\Sigma_m^+}\,$
is a local diffeomorphism,  $\,\hat P_m\,$ is necessarily 
injective on every local $\,R$-orbit of $\,\Sigma_m$. Therefore $\,iii)\,$
implies that  $\,\hat P_m|_{\Sigma_m^+}\,$ is injective. Moreover,
 from $\,i)\,$ of Lemma
\ref{FIXPOINTS} it follows that $\,\Sigma_m^+\,$ is mapped to one 
of the two connected components of $\,\mathcal P \setminus {\rm fix}(\beta)$.
By $\,\alpha$-$\beta$-equivariance of $\,\hat P_m\,$ this implies that 
$\,\Sigma_m^-:= \alpha(\Sigma_m^+)\,$ is injectively mapped to 
the other connected component. Finally it is easy
to check that $\,\k\cup {\rm fix}(\beta)\,$
is injectively mapped into $\,\{\,st=1\,\}\cup {\rm fix}(\beta)$, implying $\,iv)$. 

Now note that a local $\,R$-orbit  of $\,\Sigma_m\,$ 
either meets $\,{\rm fix}(\alpha)\,$ in a unique point or
is contained  in $\,\Sigma_m^+ \cup
\Sigma_m^-$.
As noticed, $\,\hat P_m\,$ maps different elements
of $\,{\rm fix}(\alpha)\,$ into different $\,R$-orbits of $\, \mathcal P$.
Then $\,iii)\,$ and $\,\alpha$-$\beta$-equivariance of $\,\hat P_m$,
imply $\,v)$.

For $\,vi)$, one can move towards infinity 
(topologically) along local $\,R$-orbits of $\,\Sigma_m\,$
which are closed to $\,\{e\} \times \k$, apply $\,\hat P_m\,$
and check that one is moving towards infinity along 
$\,R$-orbits in $\, \mathcal P$.
Recalling that $\,e^{2u(1+m)}\left( C(x) - 2umS(x)\right) \,$
is the second component of $\,\hat P_m$, this follows by
showing that for $\,\varepsilon>0\,$ small enough and $\,uU+aH\in 
\Sigma_m$, with $\, a< 
\frac{\pi}{4}$, such that $\,2aS(x) =  \varepsilon$,  one has
$$\lim_{u \to \infty}e^{2u(1+m)}\left( C(x) - 2umS(x)\right)= \infty\,.$$
The details of this computation are omitted.
\end{proof}

\bigskip
\begin{theorem}
\label{DOMAIN<-1}
Let $\,G=SL_2(\R)\,$ endowed with a metric $\,\nu_m$, with $\,m\leq -1$.
Then  the polar map $\, P_m|_{\Omega_m}:\Omega_m
\to G^\C\,$ is injective and consequently $\,\Omega_m \,$ is $\,L$-equivariantly
biholomorphic to $\, P_m(\Omega_m)$. This domain is 
not holomorphically convex and its envelope of holomorphy is
biholomorphic to $\,G^\C$.
\end{theorem}

\nsmallskip
\begin{proof}
$\,$Injectivity follows from $\,iv)\,$ of Lemma \ref{LEVELSET}   
and Proposition \ref{INIINI}. 
Assume by contradiction that $\,\Omega_m \,$ is holomorphically
convex, i.e. that the domain $\,P_m(\Omega_m)\,$ is Stein.
Then  the categorical
quotient   $\,P_m(\Omega_m) \qq K\,$ with respect to the $\,K$-action
is Stein (see \cite{Hn}, Sect. 6.5). Note that all local $\,K^\C$-orbits are closed in
$\,P_m(\Omega_m)\,$ and  by $\,v)\,$ of Lemma
\ref{LEVELSET} the domain $\,P_m(\Omega_m)\,$ is $\,K$-orbit-convex in $\,G^\C$.
It follows that $\,P_m(\Omega_m) \qq K \,$ is biholomorphic to 
$\, \Pi(P_m(\Omega_m))$, where  
$\,\Pi:G^\C\to G^\C/K^\C\,$ is the canonical projection.
In particular  $\,\Pi(P_m(\Omega_m))$ is Stein. 

Denote by $\,``\hat P_m(\Sigma_m)/L"\,$ 
the image of $\,\hat P_m(\Sigma_m)\,$
in $\,\mathcal P/L \,$ via the canonical projection.
Consider the  commutative diagram

$$\begin{matrix} 
  P_m(\Omega_m)  \subset G^\C  \ & \buildover{F}{\longrightarrow}  &  
  \hat P_m(\Sigma_m) \subset \mathcal P  \quad \quad \cr
                                 &                          &                        \cr
      \  \ \  \Pi\ \, \downarrow \ \        &            &    \downarrow  \quad \quad \quad \ \ \cr
                                 &                         &                         \cr
 \Pi(P_m(\Omega_m) ) \subset G^\C/K^\C  & \longrightarrow   &
                                               G \backslash \Pi(P_m(\Omega_m)) \,\cong\,
											   ``\hat P_m(\Sigma_m)/L"  \,\subset\, \mathcal P/L\,,	&	  \cr 
\end{matrix} $$

\nmedskip
where $\,F:G^\C \to {\mathcal P}\,$ is our usual quotient map.
One checks (cf. the proof of Lemma \ref{LEVELSET}) that 
$\,\hat P_m(\Sigma_m)\,$ intersects all local
$\,R$-orbits of $\,\mathcal P\,$ but one, namely 
$\,R\cdot(-1,-1)$. Moreover, for 
$\,\tilde e:=\left ( \begin{matrix}  i  &  0   \cr   
  0    &  -i	  \cr \end{matrix} \right ), $ one has  $\,F(\tilde e)=(-1,-1)$.
As a consequence $\,\Pi(P_m(\Omega_m))= G^\C/K^\C \setminus G\tilde e K^\C\,$
which, 
by Lemma \ref{ENVELOPE}, is not holomorphically convex.
This gives a contradiction, showing that  $\,\Omega_m\,$ is not 
holomorphically convex.

 For the last statement, identify $\,P_m(\Omega_m)\,$ with $\,\Omega_m\,$ and note
that its envelope of holomorphy $\,\hat \Omega_m\,$
is a Stein, $\,L$-equivariant, Riemann
domain over $\,G^\C$. One has an induced Stein, $\,G$-equivariant,
Riemann domain $\,q: \hat \Omega_m \qq K \to G^\C/ K^\C$,
where $\,\hat \Omega_m \qq K\,$ denotes the $\,K$-categorical quotient 
of $\,\hat \Omega_m\,$ (see \cite{GeIa}, Sect.$\,$3). In fact, this can be
regarded as a $\,G/\{ \pm e\}$-equivariant, Riemann domain over
$\,G^\C/ K^\C$,  
since the subgroup $\,\{ \pm (e,e)\}$ of $\, L\,$ acts trivially on $\,\Omega_m$.
Then, by Theorem$\,$7.6 in \cite{GeIa} the map $\,q\,$ is injective  and 
 Corollary$\,$3.3 in \cite{GeIa} implies that $\,\hat \Omega_m\,$ is univalent.
That is, $\,\hat \Omega_m\,$ is a Stein, $\,L$-invariant domain of $\,G^\C$.
As a consequence
$\,\hat \Omega_m \qq K\,$ is biholomorphic to $\,\Pi(\hat \Omega_m)\,$ and,
by Lemma \ref{ENVELOPE}, it necessarily coincides with $\,G^\C/ K^\C$.

Choose, on a  neighborhood $\,U\,$ of $\,eK^\C\,$ in $\,G^\C/K^\C$, a Stein
local trivialization  $\,U \times \C^* \subset G^\C\,$ of the principal $\,\C^*$-bundle
$\,\Pi:G^\C \to G^\C/K^\C$. Then one has $\,\hat \Omega_m\cap
(U \times \C^*)= \{\,(z,\lambda) \in U \times \C^* \ :\ a(z) < 
|\lambda| < b(z)  \,\}$, with $\,a < b\,$ 
functions on $\,U\,$ with values in the extended real line.
Since $\,\hat \Omega_m\cap
(U \times \C^*)\,$ is Stein, the functions
$\,\log a \,$ and $\,-\log b\,$ are plurisubharmonic.
Moreover, $\,vi)\,$
of Lemma  \ref{LEVELSET} implies that $\,\log a(z) =-\log b(z)=-\infty\,$
for $\,z\,$ close to $\,eK^\C$. Thus $\,a \equiv 0\,$ and $\,b  \equiv
\infty\,$ on $\,U$. Finally a connectedness argument using 
local trivializations covering all $\,G^\C\,$ shows that this holds
for any trivialization. Hence $\,\hat \Omega_m=G^\C\,$ as claimed.
\end{proof}

\bigskip

\section{The case $\,m>-1$}

\bigskip
Let $\,G=SL_2(\R)\,$ be endowed with one of the metrics
  $\,\nu_m$.
For $\,-1<m\le 0$, similarly to the cases considered in  section $7$,
we show that the polar map $\,P_m\,$ is a biholomorphism 
from the maximal adapted complexification $\,\Omega_m\,$ 
onto an $\,L$-invariant domain of $\,G^\C$ which is not holomorphically convex.
However, 
note  that here the projection $\,\Pi(P_m(\Omega_m))\,$ of
$\,\Omega_m\,$ to $\,G^\C/K^\C\,$ misses more than one $\,G$-orbit.

Finally, in the Riemannian cases $\,m>0\,$ the  polar map $\,P_m\,$ 
is not injective and its fibers consists of at most 
two elements. The maximal complexification
$\,\Omega_m\,$ is neither holomorphically
separable, nor holomorphically convex. For all $\,m$,
its envelope of holomorphy is shown to be biholomorphic to $\,G^\C$.


\bigskip
\begin{lem}
\label{LEVELSET2}{\rm (cf. picture 3)} Let $\,m>-1$. Then
\begin{itemize}

\item[$i)$]  different local $\,R$-orbits of
$\,\Sigma_m \cap \{\, u>0\,\}\, $ are mapped by $\, \hat P_m\,$ to different 
$\,R$-orbits of $\,\mathcal P$, 

\nsmallskip
\item[$ii)$]  local $\,R$-orbits closed to $\,\{e\} \times \k\,$ are mapped 
bijectively by $\,\hat P_m\,$ to $\,R$-orbits of  $\,{\mathcal P}$, i.e.
$\,R\,$ acts globally in an $\,R$-invariant neighborhood of $\,\{e\} \times \k$,

\nsmallskip
\item[$iii)$]  for $\,-1<m\le 0\,$ one has $\, (\hat P_m|_{\Sigma_m })^{-1}({\rm
fix} \beta) \subset {\rm fix}(\alpha)\,$ and different local $\,R$-orbits of
$\,\Sigma_m\,$ are mapped by $\, \hat P_m\,$ to different 
$\,R$-orbits of $\,\mathcal P$,

\nsmallskip
\item[$iv)$] for $\,-1<m\le 0\,$ the polar map $\,\hat P_m\,$ is injective on 
$\,\Sigma_m$,

\nsmallskip
\item[$v)$]  for $\,m>0\,$ one has $\, (\hat P_m|_{\Sigma_m })^{-1}({\rm
fix}(\beta)) \not\subset {\rm fix}(\alpha)$. More precisely $\,\Sigma_m\,$ contains
 the graph of $\,\gamma|_{\k \setminus \{0\}}$,
 with $\,\gamma\,$ defined as in Lemma \ref{FIXPOINTS}.
\end{itemize}
\end{lem}

\nsmallskip
\begin{proof} $\ $  Recall that by Lemma \ref{THESAME} the slice $\,\Sigma_m\,$ is
a  proper subdomain of $\,\Sigma^*_m.$ For $\, c \le 1\,$
let $\,\ell_c:= (\hat P_m|_{\Sigma_m^*})^{-1}(\{st=c\})$.
A similar argument as in  Lemma \ref{LEVELSET} shows that different components of 
$\,\ell_c \cap \{\,(e,\,uU+aH) \in \Sigma_m^* \ : \  a>0, \ u>0\,\}\,$ are
mapped by $\,\hat P_m\,$ to different $\,R$-orbits of
$\, \mathcal P$. One also checks that local $\,R$-orbits of $\,\Sigma_m^+:=
 \{\,(e,\,uU+aH) \in \Sigma_m \ : \  a>0, \ u>0\,\}\,$
are connected components of $\,\ell_c \cap \,\Sigma_m^+$.
Then in order to prove $\,i)\,$ it is enough to show that 
for every $\,c<1\,$ and every connected component $\, O\,$ of 
$\,\ell_c \cap \{\,(e,\,uU+aH) \in \Sigma_m^* \ : \  a>0, \ u>0\,\}\,$
the locally $\, R$-invariant set $\,O\cap \Sigma_m^+\,$ is connected.

For this  recall that (cf.$\,$Prop.$\,$\ref{DOMAIN} and Rem.$\,$\ref{THESAME})
the boundary of $\,\Sigma_m^+\,$
in  the set $\,\{\,(e,\,uU+aH) \in \Sigma_m^* \ : \ u>0, \ a>0\,\}\,$ is given by  $\, y=f(x)\,$ and 
that $\,O\,$ is a connected component of 
$\, \,4a^2S^2(x)=1-c\,$,
where
$$\,x=4u^2m^2-4a^2,\quad \,y=4u^2m^2+4a^2m\quad {\rm and} 
 \quad f(x)= \frac{C(x)}{2S'(x)}\,.$$
Since $\, 4a^2=\frac{y-x}{1+m}\,$, in the coordinates $\,x,\,y\,$
such equations read as $\, y-f(x)=0 \,$ and $\,y=\frac{(1+ m)(1-c)}{S^2(x)}+x$.
Thus it is enough to note that the function
\smallskip
$$\,\frac{(1+ m)(1-c)}{S^2(x)}+x-f(x)$$
can be rewritten as 
$$\frac{(1+ m)(1-c)}{S^2(x)}-\frac{xS(x)}{C(x)-S(x)}\,,$$
hence it is strictly decreasing by Lemma \ref{FUNCTIONS}.
This proves $\,i)$.

The analogous proof as in $\,vi)\,$ of 
\ref{LEVELSET} implies $\,ii)$.

 For $\,m=0\,$ one checks directly that 
$\,(\hat P_m|_{\Sigma_m^*})^{-1}({\rm fix} (\beta)) = {\rm fix}(\alpha)\cup
\partial \Sigma_m= \{u=0\}\cup\{a=\frac{\pi}{4} \}\,$ 
(cf. Lemma \ref{FIXPOINTS}). Thus $\,iii)\,$
holds for $\,m=0$. Now let $\,-1<m < 0\,$ and let 
$\,(e,\, X)\,$ be an element of 
$\,(\hat P_m|_{\Sigma_m^*})^{-1}({\rm fix}(\beta))\setminus {\rm fix} (\alpha)$.
That is,
$\,X=uU+aH\,$ with $\,u\not=0\,$ and
\begin{equation}
 \label{FIXED}
\frac{m}{1+m}\,\frac {C(4u^2(1+m)^2)}{S(4u^2(1+m)^2)}=\frac{C(x)}{S(x)}\,,
\end{equation}

\nsmallskip
where $\,x=4u^2m^2-4a^2>-\pi^2$.
In order to show that $\,(e,\, X)\,$ does not belong to $\,\Sigma_m\,$
we need to see that $\,y>f(x)$, i.e that (cf. Lemma$\,$\ref{FUNCTIONS}) 

$$4u^2m^2+4a^2m> \frac{xC(x)}{C(x)-S(x)}.$$

\nsmallskip
Since $\,-1 < m <0$, equation (\ref{FIXED})  above implies that
 $\,\frac{C(x)}{S(x)}<0$.
Then this inequality can be written as 

$$\left(1-\frac{S(x)}{C(x)}\right)(4u^2m^2+4a^2m)>(4u^2m^2-4a^2)$$

\noindent
which becomes 
$$4a^2m + 4a^2  > \frac{S(x)}{C(x)}(4u^2m^2+4a^2m)\,.$$

\nmedskip
Note that  $\,\frac{C(x)}{S(x)}<0\,$
also implies $\, a \not=0$.
Then, by using 
equation (\ref{FIXED}), one obtains 

$$1>\frac{S(x)}{C(x)}\frac{(4u^2m^2+4a^2m)}{4a^2(1+m)}=
\frac{S(4u^2(1+m)^2)}{C(4u^2(1+m)^2)}\bigl (\frac{u^2}{a^2}m+1 \bigr )$$

\nmedskip
which is easily checked to hold, since $\,0<\frac{S(4u^2(1+m)^2)}{C(4u^2(1+m)^2)}<1\,$
and $\,(\frac{u^2}{a^2}m+1)<1$. Finally,
a similar argument as in Lemma \ref{LEVELSET} implies the
last statement in $\,iii)$. Part $\,iv)\,$ follows from $\,iii)\,$
by using the same argument as in Lemma \ref{LEVELSET}

For $\,v)\,$ first note that $\,\Sigma_m \cap \{\,u=0\,\}\,=\,\{\,(e,\,aH) \ : \ 0\le a < \tilde a\,\}$,
with $\,\tilde a \,$ uniquely defined by $\, 0<\tilde a<\pi/2\,$ and $\,\tan (2\tilde a)=
2\tilde a\frac {1+m}{m}$.
Rewrite condition $\,**)\,$ of Proposition \ref{DOMAIN} as 
$\, x+(1+m)4a^2<\frac{xC(x)}{C(x)-S(x)}$. That is

$$(1+m)4a^2<\frac{S(x)}{2S'(x)}\,.$$

\nsmallskip
Recall that $\,\frac{S(x)}{S'(x)}\,$ is strictly increasing for $\,x>-\pi^2\,$
(Lemma \ref{FUNCTIONS}).
Then for $\,(e, \,aH) \in \Sigma_m \, $ and $\,u \in \R\,$
one has 

$$(1+m)4a^2<\frac{S(-4a^2)}{2S'(-4a^2)} \le \frac{S(x)}{2S'(x)}$$

\nsmallskip
implying that $\,\{\,(e, \,uU+aH) \ : \ u \in \R, \ 0\le a < \tilde a\,\} \subset \Sigma_m$.
Since one sees that $\,(e, \,\tilde aH) = \gamma(0)$, in order to
prove that the graph of $\,\gamma |_{\k \setminus 0}\,$ is contained in
$\, \Sigma_m\,$
it is enough to check
that  $\,\gamma\,$ is strictly decreasing (increasing) for $\,u>0\,$ ($\,u<0$).
This can be done by a direct computation showing that 
$$\,\frac{\partial \Gamma}{\partial a}\bigr|_{\{\Gamma=0\}}<0 \ \ 
{\rm and} \ \  \frac{\partial \Gamma}{\partial u}\bigr|_{\{\Gamma=0\}}<0
\ \quad \quad \ 
\Bigl ( \,\frac{\partial \Gamma}{\partial a}\bigr|_{\{\Gamma=0\}}<0  \ \ 
{\rm and } \ \  \frac{\partial \Gamma}{\partial u}\bigr|_{\{\Gamma=0\}}>0\, \Bigr ).$$
\end{proof}


\bigskip
\begin{cor}
\label{INIPART}$\ $ Let  $\,m>-1$. Then  the restrictions
$\,P_m|_{L\cdot\{\Sigma_m\cap\{u\ge 0\}\}}\,$
and $\,P_m|_{L\cdot\{\Sigma_m\cap\{u\le 0\}\}}\,$ are injective.
\end{cor}

\nsmallskip
\begin{proof}
First consider the case $\,u \ge 0$. Since $\,\hat P_m\,$ is necessarily
injective on local $\,R$-orbits, 
 $\,i)\,$ of Lemma \ref{LEVELSET2} implies that $\,\hat P_m\,$
is injective on $\,\Sigma_m^+= \Sigma_m\cap\{u> 0\}$.
Moreover, 
a straightforward computation shows that $\,\hat P_m\,$ is also injective  on 
$\,\Sigma_m \cap \{ u=0 \}$.

Then, for $ \,-1 < m \le 0$, the map   $\,\hat P_m\,$
is injective on $\, \Sigma_m\cap\{u\ge 0\}\,$ 
by $\,iii)\,$ of Lemma \ref{LEVELSET2}.
Now recall
that every $\,L$-orbit of
$\,\Omega_m\,$
meets $\,\Sigma_m\,$ in a single element and  $\,P_m\,$
is injective on such an orbit by Lemma$\,$\ref{DIFFEO}. Then 
the statement follows in the case $ \,-1 < m \le 0$.

For $\,m>0\,$ note that by $\,v)\,$ of Lemma \ref{LEVELSET2}
 the set $\,(\hat P_m|_{\Sigma_m \cap\{\,u \ge 0\,\}})^{-1}({\rm fix} (\beta))\,$
consists of a one-dimensional
manifold with two connected components.
Then in these cases one can essentially argue as follows.

Assume by contradiction that 
$\,\hat P_m(e,\,X)=\hat P_m (e,\,Y)\,$ for some
$\, (e,\,X) \,$ in $\, \{u=0\} \,$ and
$\, (e,\,Y)\,$ in $\, {\rm graph}(\gamma|_{\k \setminus 0})$.
Since $\, (e,\,X) \in \partial \Sigma^+_m$, $\, (e,\,Y) \in 
\Sigma^+_m\,$ are both in $\,\Sigma_m$, where  $\,\hat P_m\,$ is a 
local diffeomorphism,
this would imply that $\,\hat P_m|_{\Sigma_m^+}\,$
is not injective, which gives a contradiction. 

Finally the case $\,u \le 0\,$ follows from $\,\alpha$-$\beta$-equivariance of
$\,\hat P_m$.
\end{proof}

\bigskip
Let $\,\Omega\,$ be a complex $\,S^1$-manifold. We need the following 
remark on limits of local orbits for the induced local $\,\C^*$-action.
For details on induced
local actions in this particular situation we refer to \cite{HnIa}.


\bigskip
\begin{prop}
\label{OLOCONV}  Let $\,\Omega\,$ be a complex $\,S^1$-manifold
and consider the induced local
$\,\C^*$-action. Assume that there exist a sequence $\,\{x_n\}\,$ of $\,\Omega\,$ and 
an element $\,X \in
Lie(S^1) \,$ such that $\,\exp (iX) \,$ acts on every $\,x_n\,$ and  the sequences
$\,\{x_n\}$, $\, \{\exp (iX) \cdot x_n\} \,$ converge to different local
$\,\C^*$-orbits. Then $\,\Omega\,$ admits no continuous
plurisubharmonic exhaustion.
In particular it is not holomorphically convex.
\end{prop}

\nsmallskip
\begin{proof} 
Assume by contradiction that a continuous, plurisubharmonic exhaustion $\,\varphi \,$
of $\,\Omega\,$
exists. 
After integration over $\,S^1$, this function can be assumed to be 
$\,S^1$-invariant. Let $\,x\,$ and $\,y\,$ be the limit points
of $\,x_n\,$ and  of $\, \exp (iX) \cdot x_n$. 
Denote by  $\,O_x\,$ and $\,O_y\,$ the local $\,\C^*$-orbits through
$\,x\,$ and $\,y$, respectively, and choose
$\,M \in \R\, $ such that $\,\varphi(x)<M  \,$ and $\,\varphi(y)<M$.
By assumption $\,O_x\cap O_y= \emptyset$, therefore
$\,O_x\,$ is given by
$\,\{ \,\exp (\lambda X) \cdot x \ : \  a< {\rm Im}\, \lambda < b \,\}$,  with
$\,-\infty  \le a <0\,$ and $\,0<b<1$. Then $\, \exp (itX) \cdot x \to \infty\,$
as $\,t \to b$, in the sense of leaving all compact subsets of $\,X$. Thus
there exists a real $\,\tilde b\,$ with $\,0<\tilde b<b\,$ such that $\,\varphi(\exp (i \tilde b X) \cdot x) >M$.
Furthermore $\, \exp (i \tilde b X) \cdot x_n \to \exp (i \tilde b X) \cdot x$,
thus for $\,n\,$ large enough $\,\varphi(\exp (i \tilde b X) \cdot x_n) >M$, while
$\,\varphi(x_n) <M\,$ and $\,\varphi(\exp (iX) \cdot x_n) <M$.

However $\,S^1$-invariance and plurisubharmonicity of $\,\varphi\,$
imply that  the function $\,t \to \varphi(\exp (it X) \cdot x_n)\,$ is convex.
This gives a contradiction and concludes the proof.
\end{proof}

\bigskip
In fact  a similar argument yields the analogous result for actions of
compact Lie groups on holomorphically separable 
complex manifolds. The remaining cases with $\,m > -1\,$ are discussed in
the following theorem.


\bigskip
\begin{theorem}
\label{DOMAIN>-1} Let $\,G=SL_2(\R)\,$ endowed with a left invariant metric
$\,\nu_m$, with $\,m > -1$. The polar map
 $\,P_m|_{\Omega_m}:\Omega_m \to G^\C\,$ is

\begin{itemize}

\smallskip
\item[$i)$]   injective for $\,-1 < m \le0\,$ and consequently $\,\Omega_m\,$
is $\,L$-equivariantly biholomorphic to $\, P_m(\Omega_m)$.
Such domain is  
not holomorphically convex,

\nsmallskip
\item[$ii)$]   not injective for $\,m>0\,$
and its fibers
have at most two elements. The maximal complexification
$\,\Omega_m\,$ is a  Riemann domain 
over $\,G^\C\,$ which is neither holomorphically separable, nor  holomorphically convex.
\end{itemize}
\noindent
In both cases the envelope of holomorphy of $\,\Omega_m\,$ is 
biholomorphic to $\,G^\C$.
\end{theorem}

\nsmallskip
\begin{proof}
$\ $From $\,iv)\,$ of
Lemma \ref{LEVELSET2}  and Proposition \ref{INIINI} it follows that $\,P_m|_{\Omega_m}\,$  is injective for $\,-1 < m \le0$.
Assume by contradiction that $\,P_m(\Omega_m)\,$ is Stein. 
Then the analogous argument as in 
Theorem \ref{DOMAIN<-1} shows that 
$\,\Pi(P_m(\Omega_m))\,$ is Stein. 
However one checks that $\, [0,\frac{\pi}{4}]H
\in \Sigma_m$ (cf. proof of Lemma \ref{LEVELSET2}). Then, by $\,G$-invariance of 
$\,\Omega_m\,$ and $\,G$-equivariance of $\,P_m\,$,  one has  
$\,G \exp(i[0,\frac{\pi}{4}]H) K^\C \,\subset \,\Pi(P_m(\Omega_m))$. As
a consequence the closure of the Akhiezer-Gindikin domain 
$\, {\Omega}_{AG}  \,$ is contained in $\,\Pi(P_m(\Omega_m))$.
Since one sees that $\,\Pi(P_m(\Omega_m))\, \not =\, G^\C/K^\C\,$
(cf. proof of Lemma \ref{LEVELSET2}), Lemma \ref{ENVELOPE} 
implies that $\,\Pi(P_m(\Omega_m))\,$ is not Stein.
This gives a contradiction, implying $\,i)$.

For $\,ii)\, $ note that by $\,v)\,$ of Lemma \ref{LEVELSET2}
there exists an element $\,(e,\,X)\,$ in $\,\Sigma_m \cap(\hat P_m)^{-1}({\rm fix \beta})
\,$ with  $\,X_\k \not=0$.
Then by $\,\alpha $-$ \beta $-equivariance of $\,\hat P_m\,$ one has
$\,\hat P_m(e,\,X_\k+X_\a)=\hat P_m(e,\,-X_\k+X_\a)$,
showing that $\,\hat P_m|_{\Sigma_m}\,$ is not injective. Then,
 by Proposition 
\ref{INIINI} the polar map  $\,P_m|_{\Omega_m}\,$  is not injective as well.

Since $\,\Sigma_m=(\Sigma_m\cap\{u\ge 0\}) \cup (\Sigma_m\cap\{u\le 0\})$,
from Corollary \ref{INIPART} and Proposition \ref{INIINI}
it follows that the fibers of $\,P_m\,$
consist either of one point or two points $\,p^+\,$ and $\,p^-$.
In the second case, necessarily $\,p^+ \in L \cdot\Sigma_m^+\,$ and $\,p^- \in
L \cdot\Sigma_m^-$,
where  $\,\Sigma_m^-= \alpha(\Sigma_m^+)=\{\, (e,\,uU+aH) \,:\, u<0, \ a>0\,\}$.

Now let $\,(e,\,[0,\tilde a)H) =  \Sigma_m \cap \a^+\,$ and
note that local $\,\C^*$-orbits through $\,(e,\,tH)\,$ accumulate,
for $\,t \to \tilde a$, to different local $\,\C^*$-orbits.
Indeed this holds true for their images in the quotient $\, \Sigma_m \cong \Omega/L$.
Such images are
given by the $\,R$-orbits through $\,(e,\,tH)\,$ for which, chosen $\,\varepsilon\,$
small enough,
$\,(\{e\} \times \exp (i \varepsilon U)) \cdot (e,\,tH)\,$ and 
$\,(\{e\} \times \exp (-i \varepsilon U)) \cdot (e,\,tH)\,$ accumulate to
different local $\,R$-orbits. Namely the two connected components
of $\,\ell_{\tilde c} \cap \Sigma_m$, with $\,\tilde c\,$ such that
$\, \hat P_m(e,\,\tilde aH) \in \{\,st=\tilde c\,\}\,$ (cf. the proof of Lemma
\ref{LEVELSET2}). Then Proposition \ref{OLOCONV} implies that 
$\,\Omega_m\,$ is not holomorphically convex.

Assume by contradiction that 
$\,\Omega_m\,$ is holomorphically separable. Then it
embeds in its envelope of holomorphy $\,\hat \Omega_m$,
which is a non-univalent, Stein, $\,L$-equivariant, Riemann domain over $\,G^\C$.
By Corollary$\,$3.3 in \cite{GeIa}, the induced 
Stein, $\,G$-equivariant, Riemann domain $q:\hat \Omega_m \qq K  \to G^\C/K^C\,$
is also non-univalent. However the subgroup $\,\{\pm(e,e)\}\,$
of $\,L\,$ acts trivially on $\,\hat \Omega_m$, thus 
$\,\hat \Omega_m \qq K\,$ can be regarded as a 
 $\,G/\{\pm e\}$-equivariant, Riemann domain. Then Theorem$\,$7.6
in \cite{GeIa} implies that $\,q\,$ is injective, giving a contradiction.

Finally, in view of $\,ii)\,$ of Lemma \ref{LEVELSET2},
the last statement follows from the same argument as in Theorem \ref{DOMAIN<-1}.
\end{proof}

\bigskip

\section{ Appendix}

\bigskip
\noindent

Here we carry out the computation used in  section 4.
For this note  that  $\, [U,H]\,=\,2W\,$, $\,[U,W]\,=\,-2H\,$ and $\,[H,W]\,=\,-2U$.
Thus
\smallskip
\begin{align*}
\ad(i(umU-aH))\,U\,= & \,[i(umU-aH), U]\,=\,2aiW, \\  
\ad^2(i(umU-aH))\,U\,=&\,[i(umU-aH),2aiW ]\,=\\
&\,-4a(aU-umH),\\
\ad^3(i(umU-aH))\,U\,=&\,[i(umU-aH),-4a(aU-umH) ]\,= \\
& \, 8a(u^2m^2-a^2)iW,\\
\ad^4(i(umU-aH))\,U\,=&\,[i(umU-aH),8a(u^2m^2-a^2)iW ]\,=\\
& \, -16a(u^2m^2-a^2)(aU-umH).
\end{align*} 
         
\noindent
Then, by recalling that $\,\Ad_{\exp}=e^{\ad}$, one has 

$$\Ad_{\exp i(umU-aH)}U \,= e^{\ad(i(umU-aH))}U = U\,+\,2aiW \,
- \,\frac{4a}{2!}\,(aU-umH) \,+$$

$$ \frac{8a(u^2m^2-a^2)}{3!}\,iW
\,- \, \frac{16a(u^2m^2-a^2)}{4!}\,(aU-umH) \,+ \,\dots \,= $$

$$U \,-\, 4a \frac{\cosh \sqrt{x} -1}{x}\,(aU-umH)\,+ \,  2a
\frac{\sinh \sqrt{x}}{\sqrt{x}}\,iW  =$$

$$\left ( 1-4a^2\frac{C(x)-1}{x} \right) U \ + \ 
 4aum\frac{C(x)-1}{x} \, H\  +\ 
2aS(x) \, iW\,,$$

\medskip
\noindent
where $\,x= 4u^2m^2 -4a^2$.
Similarly, for the second vector one has
\smallskip
\begin{align*}
\ad(i(umU-aH))H\,=&\,[i(umU-aH), H]\,=\,2umiW,\\
\ad^2(i(umU-aH))H\,=&\,[i(umU-aH),2umiW ]\,=\\
&\,-4um(aU-umH),\\
\ad^3(i(umU-aH))H\,=&\,[i(umU-aH),-4um(aU-umH)]\,= \\
&\, 8um(u^2m^2-a^2)iW,\\
\ad^4(i(umU-aH))H\, =&\,[i(umU-aH),8um(u^2m^2-a^2)iW]\,=\\
&\,-16um(u^2m^2-a^2)(aU-umH).
\end{align*}

\medskip
\noindent
Therefore $\,\Ad_{\exp i(umU-aH)}H \,= \,e^{\ad(i(umU-aH))}H \, = \,$
\smallskip
$$ H \, + \, 2umiW
 \, -\, \frac{4um}{2!}\,(aU-umH) \,+\,$$
\smallskip
$$ \frac{8um(u^2m^2-a^2)}{3!}\,iW \, -\,
\frac{16um(u^2m^2-a^2)}{4!}\,(aU-umH) \, + \, \dots \,=$$
\smallskip
$$H \,-\, 4um \frac{\cosh \sqrt{x} -1}{x}\,(aU-umH) \,+ \, 2um
\frac{\sinh \sqrt{x}}{\sqrt{x}}\,iW\,
 =$$
\smallskip
$$-4aum\frac{C(x)-1}{x} \, U \ +\  
\left ( 1+4u^2m^2\frac{C(x)-1}{x} \right) H\  + \
2umS(x) \,iW\,.$$

\nmedskip
For the third vector one has
\medskip
\begin{align*}
\ad(i(umU-aH))\,W\,= &\,[i(umU-aH), W]\,=\,2i(aU-umH), \\
\ad^2(i(umU-aH))\,W\,= & \,[i(umU-aH),2i(aU-umH)]\,=\\
&\,4(u^2m^2-a^2)\,W, \\
\ad^3(i(umU-aH))\,W\,= &\,[i(umU-aH),4(u^2m^2-a^2)W] \,= \\
&\,8(u^2m^2-a^2)i(aU-umH), \\
\ad^4(i(umU-aH))\,W\,= &\,[i(umU-aH),8(u^2m^2-a^2)i(aU-umH)]\,= \\
&\,16(u^2m^2-a^2)^2W.
\end{align*}

\medskip
\noindent
Therefore $\,\Ad_{\exp i(umU-aH)}W \,=\, e^{\ad(i(umU-aH))}W \,=\,$ 
\smallskip
$$W \,+ \,2i(aU-umH) \,+\,
 \frac{4(u^2m^2-a^2)}{2!}\,W \,+$$
\smallskip
$$\frac{8(u^2m^2-a^2)}{3!}\,i(aU-umH) \,+ \,
\frac{16(u^2m^2-a^2)^2}{4!}\,W \,+ \,\dots \,= $$
\smallskip
$$2  \frac{\sinh\sqrt{x}}{\sqrt{x}} \,i(aU - 2umH) \, +\, 
\cosh \sqrt{x}\, W\,=$$
\smallskip
$$2aS(x) \,iU \ -\ 2umS(x) \,iH \ + \ C(x)\, W.$$

\nmedskip
For the fourth vector one has
\medskip
\begin{align*}
\ad(-i(umU-aH))(-imU)\,=&\,[-i(umU-aH), -imU]\,=\,-2amW,\\
\ad^2(-i(umU-aH))(-imU)\,=&\,[-i(umU-aH),-2amW ]\,=\\
&\,4ami(aU-umH),\\
\ad^3(-i(umU-aH))(-imU)\,=&\,[-i(umU-aH),4ami(aU-umH) ]\,= \\
&\,-8am(u^2m^2-a^2)W,\\
\ad^4(-i(umU-aH))(-imU)\,=&\,[-i(umU-aH),-8am(u^2m^2-a^2)W]\,= \\
&\,16am(u^2m^2-a^2)i(aU-umH).
\end{align*}

\medskip
\noindent
Therefore $\, \sum_{l=0}^\infty \frac{(-1)^l}{(l+1)!}\ad^l(i(-umU+aH))(-imU) \ +  
i(1+m)U\, =$
$$ iU \, +\, \frac{2am}{2!}\,W \, + \,\frac{4am}{3!}\,i(aU-umH) \, + \,$$

$$ \frac{8am(u^2m^2-a^2)}{4!}\,W
\, +\, \frac{16am(u^2m^2-a^2)}{5!}\,i(aU-umH) \, +\, \dots \,= $$

$$iU \,+\, 4am \frac{\sinh \sqrt{x}/\sqrt{x}-1}{x}\,i(aU-umH) \,+\,
2am \frac{\cosh \sqrt{x}-1}{x}\,W  =$$

$$\left ( 1+4a^2m\frac{S(x)-1}{x} \right) \,iU\  - \ 4aum^2\frac{S(x)-1}{x}\,iH \ + \ 
2am\frac{C(x)-1}{x} \,W\,.$$

\nmedskip
For the fifth vector one has
\begin{align*}
\ad(-i(umU-aH))\, iH\,=&\,[-i(umU-aH), iH]\,=\,2umW, \\
\ad^2(-i(umU-aH))\,iH\,=&\,[-i(umU-aH),2umW]\,=\\
&\,-4umi(aU-umH),\\
\ad^3(-i(umU-aH))\,iH\,=&\,[-i(umU-aH),-4umi(aU-umH)]\,= \\
&\,8um(u^2m^2-a^2)W,\\
\ad^4(-i(umU-aH))\,iH\,=&\,[-i(umU-aH),8um(u^2m^2-a^2)W]\,=\\
&\,-16um(u^2m^2-a^2)i(aU-umH).
\end{align*}
\medskip
\noindent
Therefore $\,\sum_{l=0}^\infty \frac{(-1)^l}{(l+1)!}\ad^l(i(-umU+aH))(iH)\,=\,$
$$ 
iH \,- \,\frac{2um}{2!}\,W \,-\, \frac{4um}{3!}\,i(aU - umH) \,- $$
\smallskip
$$- \,\frac{8um(u^2m^2-a^2)}{4!}\,W
\,-\, \frac{16um(u^2m^2-a^2)}{5!}\,i (aU-umH) \,-\, \dots \,= $$
\smallskip
$$iH\, -\,4um \frac{\sinh \sqrt{x}/\sqrt{x}-1}{x}\, i(umH-aU)
\,-\, 2um \frac{\cosh \sqrt{x} -1}{x}\, W  =$$
\smallskip
$$-4aum \frac{S(x)-1}{x}\, iU \ +  \ 
\left ( 1+4u^2m^2\frac{S(x)-1}{x} \right) \, \,iH \ -\  2um\frac{C(x)-1}{x} \,W\,,$$

\nmedskip
For the sixth vector one has
\nmedskip
\begin{align*}
\ad(-i(umU-aH))\,iW\,=& \,[-i(umU-aH),iW]\,=\,2(aU-umH),\\
\ad^2(-i(umU-aH))\,iW\,= &\,[-i(umU-aH), 2(aU-umH)]\,=\\
&\,4(u^2m^2-a^2)iW,
\end{align*}
\begin{align*}
\ad^3(-i(umU-aH))\,iW\,=&\,[-i(umU-aH),4(u^2m^2-a^2)iW]\,= \\
&\,8(u^2m^2-a^2)(aU-umH),\\
\ad^4(-i(umU-aH))\,iW\,=&\,[-i(umU-aH),8(u^2m^2-a^2)(aU-umH)]\,=\\
&\,16(u^2m^2-a^2)^2iW.
\end{align*}
\noindent
Therefore $\,\sum_{l=0}^\infty \frac{(-1)^l}{(l+1)!}\ad^l(i(-umU+aH))(iW)\,=$
\nmedskip
$$ iW \,-\,
\frac{2}{2!}\,(aU-umH) \,+\, \frac{4(u^2m^2-a^2)}{3!}\, iW \,+$$
\smallskip
$$ - \frac{8(u^2m^2-a^2)}{4!}\, (aU-umH)
\, + \, \frac{16(u^2m^2-a^2)^2}{5!}\,iW \, -\, \dots \,= $$
\smallskip
$$-\,2\frac{\cosh\sqrt{x}-1}{x} \,(aU -umH)
\,+\,  \frac{\sinh \sqrt{x}}{\sqrt{x}} \,iW\,=$$
\smallskip
$$ \,-2a\frac{C(x)-1}{x} \,U \ + \ 
2um\frac{C(x)-1}{x}\,H\ +\  S(x)\,iW.$$

\bigskip


\end{document}